\RequirePackage{fix-cm}
\documentclass{svjour3}                     
\smartqed  
\usepackage{graphicx,verbatim,xcolor,soul,relsize}
\usepackage{verbatim,xcolor,soul,bm,stmaryrd}
\usepackage{amsmath,lineno,amssymb,mathtools,url}

\newtheorem{dfn}{Definition}
\newtheorem{pro}[dfn]{Problem}
\newtheorem{thm}[dfn]{Theorem}
\newtheorem{exa}[dfn]{Example}

\newtheorem{cor}[dfn]{Corollary}
\newtheorem{lem}[dfn]{Lemma}
\newtheorem{prop}[dfn]{Proposition}
\counterwithin{dfn}{section}

\def\pbc{{\color[rgb]{0,0.5,1}{p_{23}}}}
\def\pac{{\color[rgb]{1,0.5,0}{p_{13}}}}
\def\pab{{\color[rgb]{0.5,0,0.5}{p_{12}}}}
\def\pda{{\color[rgb]{0,0,1}{p_{01}}}}
\def\pdb{{\color[rgb]{0,0.5,0}{p_{02}}}}
\def\pdc{{\color[rgb]{1,0,0}{p_{03}}}}

\def\rac{{\color[rgb]{1,0.5,0}{r_{13}}}}
\def\rab{{\color[rgb]{0.5,0,0.5}{r_{12}}}}
\def\rda{{\color[rgb]{0,0,1}{r_{01}}}}
\def\rdb{{\color[rgb]{0,0.5,0}{r_{02}}}}
\def\rdc{{\color[rgb]{1,0,0}{r_{03}}}}

\newcommand{\R}{\mathbb R}
\newcommand{\Z}{\mathbb Z}

\newcommand{\de}{\delta}

\newcommand{\La}{\Lambda}
\newcommand{\si}{\sigma}
\newcommand{\ep}{\varepsilon}

\newcommand{\ti}{\tilde}

\newcommand{\inv}{\mathrm{Inv}}

\newcommand{\VF}{\mathrm{VF}}
\newcommand{\CF}{\mathrm{CF}}

\newcommand{\RF}{\mathrm{RF}}
\newcommand{\RI}{\mathrm{RI}}

\newcommand{\RIS}{\mathrm{RIS}}

\newcommand{\LIS}{\mathrm{LIS}}

\newcommand{\ws}{\hfill $\square$}
\newcommand{\bt}{\hfill $\blacktriangle$}
\newcommand{\bs}{\hfill $\blacksquare$}
\newcommand{\lra}{\leftrightarrow}

\newcommand{\vl}{\,:\,}
\newcommand{\matfour}[4]{\left(\begin{array}{ccc}
#1 & #2 \\ #3 & #4 \end{array}\right)}
\newcommand{\mat}[6]{\left(\begin{array}{ccc}
#1 & #2 & #3 \\ #4 & #5 & #6 \end{array}\right)}

\begin{document}

\title{A complete isometry classification of 3-dimensional lattices} 


\author{Vitaliy Kurlin 
}


\institute{V.Kurlin \at 
Computer Science, University of Liverpool, UK \email{vitaliy.kurlin@liverpool.ac.uk}
}

\date{Received: date / Accepted: date}
\maketitle

\begin{abstract}
A periodic lattice in Euclidean 3-space is the infinite set of all integer linear combinations of basis vectors.
Any lattice can be generated by infinitely many different bases.
This ambiguity was only partially resolved, but standard reductions remained discontinuous under perturbations modelling crystal vibrations.
This paper completes a continuous classification of 3-dimensional lattices up to Euclidean isometry (or congruence) and similarity (with uniform scaling).
\medskip

The new homogeneous invariants are uniquely ordered square roots of scalar products of four vectors whose sum is zero and all pairwise angles are non-acute.
These root invariants continuously change under perturbations of basis vectors.
The geometric methods extend the past work of Delone, Conway and Sloane.


\end{abstract}

\section{The hard problem to continuously classify lattices up to isometry}
\label{sec:intro}

We extend the continuous isometry classification of 2-dimensional lattices \cite{kurlin2022mathematics} to dimension 3.
A \emph{lattice} $\La\subset\R^n$ consists of integer linear combinations of basis vectors $v_1,\dots,v_n$.
This basis spans a parallelepiped called a \emph{unit cell} $U\subset\R^n$.
\medskip

The problem to classify lattices up to isometry is motivated by periodic crystals whose structures are determined in a rigid form.
Hence the most natural equivalence of crystals is rigid motion.
We start from general isometries that also include mirror reflections because the sign of a lattice similar to \cite[Definition~3.4]{kurlin2022mathematics} easily distinguishes mirror images.
As in $\R^2$, the space of lattices up to rigid motion in $\R^3$ is a 2-fold cover of the smaller Lattice Isometry Space $\LIS(\R^3)$.  
\medskip

The previous work \cite[section~1]{kurlin2022mathematics} provided important motivations for a continuous classification problem, which we state below for 3-dimensional lattices.

\begin{pro}[continuous classification of 3D lattices]
\label{pro:map}
Find an invariant $I:\LIS(\R^3)\to\inv$ mapping the Lattice Isometry Space to a simpler space such that
\smallskip

\noindent
(\ref{pro:map}a) 
\emph{invariance} : $I(\La)$ is independent of a lattice basis and is preserved under isometry of $\R^3$, so $I$ has no false negatives : if $\La\cong\La'$ then $I(\La)=I(\La')$;
\smallskip

\noindent
(\ref{pro:map}b) 
\emph{completeness} : if $I(\La)=I(\La')$, then $\La,\La'$ are isometric, so $I$ has no false positives and defines a bijection (or a 1-1 map) $I:\LIS\to \inv=I(\LIS)$;
\smallskip

\noindent
(\ref{pro:map}c)
\emph{continuity} : $I(\La)$ is continuous under perturbations of a basis of $\La$;
\smallskip

\noindent
(\ref{pro:map}d)
\emph{computability} : 
$I(\La)$ can be explicitly computed from a suitable basis of $\La$;
\smallskip

\noindent
(\ref{pro:map}e)
\emph{inverse design} : a basis of $\La$ can be explicitly reconstructed from $I(\La)$.
\bs 
\end{pro}

About 30 years ago John Conway and Neil Sloane published a series of seven papers on low-dimensional lattices.
The most relevant for Problem~\ref{pro:map} is \cite[item 1 on page 55]{conway1992low} saying that certain lattice invariants (conorms) `vary continuously with the lattice'.
Unfortunately, there was no further discussion of continuity and even no rigorous statement of invariance, because the above invariants should be considered up to different permutations depending on the Voronoi type of a lattice, see Lemmas~\ref{lem:V1superbases}-\ref{lem:V5superbases}.
Otherwise Problem~\ref{pro:map} might have been solved in 1992.
\medskip

\begin{figure}[h]
\caption{
Vectors of an obtuse superbase of a lattice $\La\subset\R^3$ have ordered scalar products that form the root invariants continuously parameterising the Lattise Isometry Space $\LIS(\R^3)$.}
\label{fig:lattice_classification}
\medskip

\includegraphics[width=1.0\textwidth]{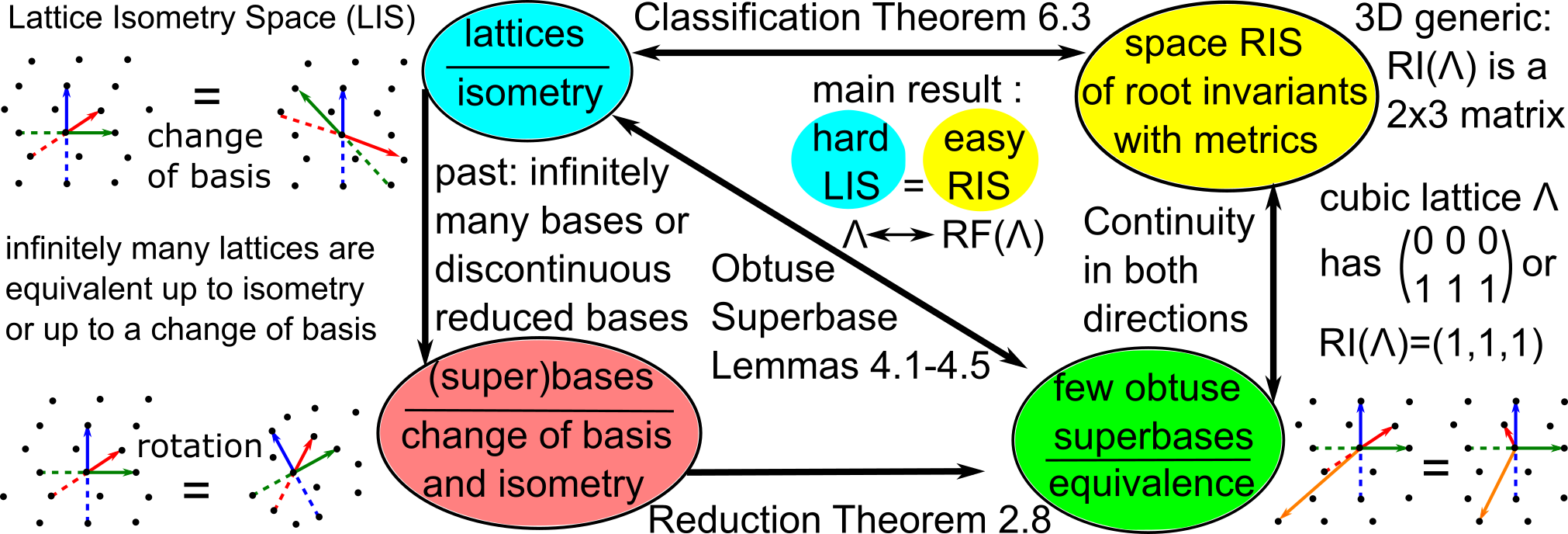}
\end{figure}

In $\R^2$, \cite[Problem~1.1]{kurlin2022mathematics} was stated and solved for stronger conditions (1.1c)-(1.1d) requiring a continuous and computable metric on lattices.
This metric part of Problem~\ref{pro:map} is postponed to the next paper, because the invariant part is already hard in $\R^3$.
The orientation-aware equivalences (rigid motion and orientation-preserving similarity) are also postponed for future work.
Fig.~\ref{fig:lattice_classification} summarises the past obstacles and a full solution to Problem~\ref{pro:map}.
The space $\inv$ will be the root invariant space ($\RIS$) of root invariants consisting of up to six parameters.

\vspace*{-2mm}
 
\section{Main definitions and an overview of past work and new results}
\label{sec:review}

The previous work defined the main concepts for any dimension $n\geq 2$ in \cite[section~2]{kurlin2022mathematics}.
For simplicity, we remind these concepts only for $n=3$.
Any point $p$ in Euclidean space $\R^3$ can be represented by the vector from the origin $0\in\R^n$ to $p$.
This vector is also denoted by $p$,
An equal vector $p$ can be drawn at any initial point.
The \emph{Euclidean} distance between points $p,q\in\R^3$ is $|p-q|$.

\begin{dfn}[a lattice $\La$, a primitive unit cell $U$]
\label{dfn:lattice_cell}
Let vectors $v_1,v_2,v_3$ form a linear {\em basis} in $\R^3$ so that 
any vector $v\in\R^3$ can be written as $v=c_1 v_1+c_2 v_2+c_3 v_3$ for some real $c_i\in\R$, and if $v=0$ then $c_1=c_2=c_3=0$.
A {\em lattice} $\La$ in $\R^3$ consists of 
$c_1 v_1+c_2 v_2+c_3$ for $c_i\in\Z$.
The parallelepiped $U(v_1,v_2,v_3)=\left\{ c_1 v_1+c_2 v_2+c_3 v_3 \vl c_i\in[0,1) \right\}$ is a \emph{primitive unit cell} of $\La$.
\bs
\end{dfn}

The conditions $0\leq c_i<1$ on the coefficients $c_i$ above guarantee that the copies of unit cells $U(v_1,v_2,v_3)$ translated by all $v\in\La$ are disjoint and cover $\R^3$.

\begin{dfn}[orientation, isometry, rigid motion, similarity]
\label{dfn:isometry}
For a basis $v_1,v_2,v_3$ of $\R^3$, the \emph{signed volume} of $U(v_1,v_2,v_3)$ is the determinant of the $3\times 3$ matrix with columns $v_1,v_2,v_3$.
The sign of this $\det(v_1,v_2,v_3)$ can be called an \emph{orientation} of the basis $v_1,v_2,v_3$. 
An \emph{isometry} is any map $f:\R^3\to\R^3$ such that $|f(p)-f(q)|=|p-q|$ for any $p,q\in\R^3$.
The unit cells $U(v_1,v_2,v_3)$ and $U(f(v_1),f(v_2),f(v_3))$ have non-zero volumes with equal absolute values.
If these volumes have equal signs, $f$ is \emph{orientation-preserving}, otherwise $f$ is \emph{orientation-reversing}.
Any orientation-preserving isometry $f$ is a composition of translations and rotations, and can be included into a continuous family of isometries $f_t$ (a \emph{rigid motion}), where $t\in[0,1]$, $f_0$ is the identity map and $f_1=f$.
A \emph{similarity} is a composition of isometry and uniform scaling $v\mapsto sv$ for a fixed scalar $s>0$. 
\bs
\end{dfn}

Any lattice $\La$ can be generated by infinitely many bases or unit cells.
This ambiguity was traditionally resolved by a reduced basis, which can be defined in several ways \cite{gruber1989reduced}.
All these reduced bases including Niggli's basis \cite{niggli1928krystallographische} are discontinuous under perturbations, which was highlighted in \cite[section~1]{edels2021}, see an example extendable to any dimension by adding long orthogonal basis vectors in \cite[Fig.~3]{kurlin2022mathematics} and a formal proof in \cite[Theorem~15]{widdowson2022average}.
Experimentally, discontinuity of Niggli's basis was demonstrated in the seminal work \cite{andrews1980perturbation} and motivated the subsequent progress of Larry Andrews and Herbert Bernstein \cite{andrews1988lattices,andrews2014geometry,mcgill2014geometry,andrews2019selling} in Problem~\ref{pro:map}.
\medskip

The proposed solution is based on the \emph{Voronoi domain} \cite{voronoi1908nouvelles}, also called the \emph{Wigner-Seitz cell}, \emph{Brillouin zone} or \emph{Dirichlet cell}.
We use the word \emph{domain} to avoid a confusion with a unit cell in Definition~\ref{dfn:lattice_cell}.
Though the Voronoi domain can be defined for any point of a lattice, it suffices to consider only the origin $0$.

\begin{dfn}[Voronoi domain $V(\La)$]
\label{dfn:Voronoi_vectors}
The \emph{Voronoi domain} of a lattice $\La\subset\R^3$ is the neighbourhood $V(\La)=\{p\in\R^3: |p|\leq|p-v| \text{ for any }v\in\La\}$ of the origin $0\in\La$ consisting of all points $p$ that are non-strictly closer to $0$ than to other points $v\in\La$.
A vector $v\in\La$ is a \emph{Voronoi vector} if the bisector hyperspace $H(0,v)=\{p\in\R^n \vl p\cdot v=\frac{1}{2}v^2\}$ between 0 and $v$ intersects $V(\La)$.
If $V(\La)\cap H(0,v)$ is a 2-dimensional face of $V(\La)$, then $v$ is called a \emph{strict} Voronoi vector. 
\bs
\end{dfn}

Voronoi \cite{voronoi1908nouvelles} proved any lattice $\La\subset\R^3$ has one of the Voronoi types below:
\smallskip

\noindent
Voronoi type $V_1$: 
a truncated octahedron;
\smallskip

\noindent
Voronoi type $V_2$: 
a hexa-rhombic dodecahedron;
\smallskip

\noindent
Voronoi type $V_3$: 
a rhombic dodecahedron;
\smallskip

\noindent
Voronoi type $V_4$: 
a hexagonal prism;
\smallskip

\noindent
Voronoi type $V_5$: 
a cuboid (an orthogonal parallelepiped or a rectangular box).
\medskip

Any lattice is determined by its Voronoi domain by \cite[Lemma~A.1]{kurlin2022mathematics}.
However, the combinatorial structure of $V(\La)$ is discontinuous under perturbations.
Almost any perturbation of an orthogonal basis in $\R^3$ (whose lattice has a cuboid Voronoi domain) gives a generic lattice whose Voronoi domain of type $V_1$.
Hence any integer-valued descriptors of $V(\La)$ such as the numbers of vertices or edges are always discontinuous and unsuitable for continuous quantification of similarities between arbitrary crystals or periodic point sets.
\medskip

Optimal geometric matching of Voronoi domains with a shared centre led \cite{mosca2020voronoi} to two continuous metrics (up to orientation-preserving isometry and similarity) on lattices.
The minimisation over infinitely many rotations was implemented in \cite{mosca2020voronoi} by sampling and gave approximate algorithms for these metrics.
The complete invariant isoset \cite{anosova2021isometry} for periodic point sets in $\R^n$ has a continuous metric that can be approximated \cite{anosova2021introduction} with a factor $O(n)$.  
The metric on invariant density functions \cite{edels2021} required a minimisation over $\R$, so far without approximation guarantees.
\medskip


Lemma~\ref{lem:Voronoi_vectors} shows how to find all Voronoi vectors of any lattice $\La\subset\R^n$. 
The doubled lattice is $2\La=\{2v \vl v\in\La\}$.
Vectors $u,v\in\La$ are called \emph{$2\La$-equivalent} if $u-v\in 2\La$.
Then any vector $v\in\La$ generates its $2\La$-class $v+2\La=\{v+2u \vl u\in\La\}$, which is $2\La$ translated by $v$ and containing $-v$.
All classes of $2\La$-equivalent vectors form the quotient space $\La/2\La$.
 
\begin{lem}[a criterion for Voronoi vectors {\cite[Theorem~2]{conway1992low}}]
\label{lem:Voronoi_vectors}
For any lattice $\La\subset\R^n$, a non-zero vector $v\in\La$ is a Voronoi vector of $\La$ if and only if $v$ is a shortest vector in its $2\La$-class $v+2\La$.
Also, $v$ is a strict Voronoi vector if and only if $\pm v$ are the only shortest vectors in the $2\La$-class $v+2\La$.
\bt
\end{lem}

We use the notations from \cite{conway1992low}, though obtuse superbases and their conorms were studied earlier by Selling \cite{selling1874ueber} for $n=3$ and Delone for any $n\geq 2$ \cite{delone1937geometry}.

\begin{dfn}[obtuse superbase, conorms $p_{ij}$]
\label{dfn:conorms}
For any basis $v_1,v_2,v_3$ in $\R^n$, the \emph{superbase} includes the vector $v_0=-v_1-v_2-v_3$.
The \emph{conorms} $p_{ij}=-v_i\cdot v_j$ are the negative scalar products of the vectors above. 
The superbase is \emph{obtuse} if all conorms $p_{ij}\geq 0$, so all angles between vectors $v_i,v_j$ are non-acute for distinct indices $i,j\in\{0,1,2,3\}$.
The superbase is \emph{strict} if all $p_{ij}>0$.
\bs
\end{dfn}

\cite[formula (1)]{conway1992low} has a typo initially defining $p_{ij}$ as exact Selling parameters, but later Theorems 3,~7,~8 use the non-negative conorms $p_{ij}=-v_i\cdot v_j\geq 0$.
\medskip

The indices of a conorm $p_{ij}$ are distinct and unordered.
We set $p_{ij}=p_{ji}$ for all $i,j$.
Any superbase of $\R^3$ has six conorms $p_{12},p_{13},p_{23},p_{01},p_{02},p_{03}$.

\begin{dfn}[partial sums $v_S$, vonorms $v_S^2$]
\label{dfn:vonorms}
Let a lattice $\La\subset\R^n$ have a superbase $B=\{v_0,v_1,v_2,v_3\}$. 
For any proper subset $S\subset\{0,1,2,3\}$ of indices,
 consider its complement $\bar S=\{0,1,2,3\}-S$ and the \emph{partial sum} $v_S=\sum\limits_{i\in S} v_i$ whose squared lengths $v_S^2$ are called the \emph{vonorms} of $B$ and can be expressed as $$v_S^2=(\sum\limits_{i\in S} v_i)(-\sum\limits_{j\in\bar S}v_j)=-\sum\limits_{i\in S,j\in\bar S}v_{j}\cdot v_j=\sum\limits_{i\in S,j\in\bar S}p_{ij}.\leqno{(\ref{dfn:vonorms}a)}$$
$$\text{For example },
v_i^2=p_{ij}+p_{ik}+p_{il} \text{ for any unordered triple } \{j,k,l\}=\{0,1,2,3\}-\{i\},$$ 
$$v_{ij}^2=(v_i+v_j)^2=(-v_k-v_l)^2=p_{ik}+p_{il}+p_{jk}+p_{jl} \text{ for }\{k,l\}=\{0,1,2,3\}-\{i,j\}.$$
For instance, $v_0^2=p_{01}+p_{02}+p_{03}$.
The six conorms are conversely expressed as 
$$p_{ij}=\dfrac{1}{2}(v_i^2+v_j^2-v_{ij}^2) \text{ for any distinct indices } i,j\in\{0,1,2,3\}.\leqno{(\ref{dfn:vonorms}b)}$$
The seven vonorms above have the relation
$v_0^2+v_1^2+v_2^2+v_3^2=v_{01}^2+v_{02}^2+v_{03}^2$.
\bs
\end{dfn}

Lemma~\ref{lem:partial_sums} will help 
classify obtuse superbases for all five Voronoi domains.

\begin{lem}[Voronoi vectors $v_S$ {\cite[Theorem~3]{conway1992low}}]
\label{lem:partial_sums}
For any obtuse superbase $v_0,v_1,v_2,v_3$ of a lattice, all partial sums $v_S$ from Definition~\ref{dfn:vonorms} split into seven symmetric pairs $v_S=-v_{\bar S}$, which are Voronoi vectors representing distinct $2\La$-classes in $\La/2\La$.
All Voronoi vectors $v_S$ are strict if and only if all $p_{ij}>0$.
\bt
\end{lem}

By Conway and Sloane \cite[section~2]{conway1992low}, any lattice $\La\subset\R^n$ that has an obtuse superbase is called a \emph{lattice of Voronoi's first kind}.
It turns out that any lattice in dimensions 2 and 3 is of Voronoi's first kind by Theorem~\ref{thm:reduction}, likely for any $n\geq 4$ because higher dimensions have `more space' for obtuse superbases.

\begin{thm}[reduction to an obtuse superbase]
\label{thm:reduction}
Any lattice $\La\subset\R^3$ has an obtuse superbase $\{v_0,v_1,v_2,v_3\}$ so that  all conorms $p_{ij}=-v_i\cdot v_j\geq 0$. 
\bt
\end{thm}

Conway and Sloane in \cite[section~7]{conway1992low} attempted to prove Theorem~\ref{thm:reduction} for $n=3$ by example whose details are corrected after the updated proof in appendix~\ref{sec:proofs}.

\section{Voforms and coforms of an obtuse superbase of a 3D lattice}
\label{sec:forms3d}

\begin{figure}[h]
\caption{\textbf{Left}: the Fano plane is a set of seven nodes arranged in triples shown by six lines and one circle.
\textbf{Middle}: nodes of the voform $\VF(\La)$ are marked by vonorms $v_i^2$ and $v_{ij}^2$.
\textbf{Right}: nodes of the coform $\CF(\La)$ are marked by conorms $p_{ij}$ and 0.}
\medskip

\includegraphics[height=27mm]{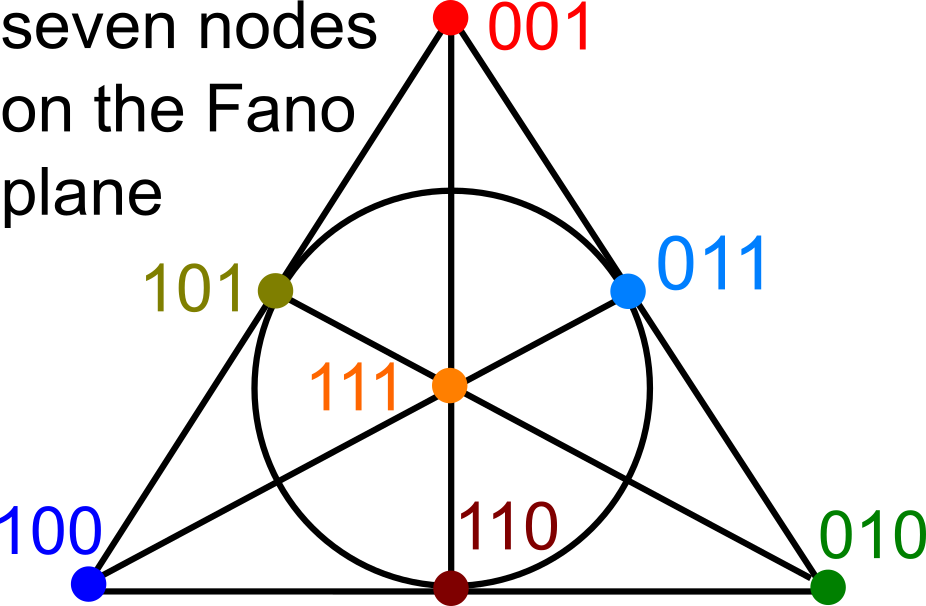}
\includegraphics[height=27mm]{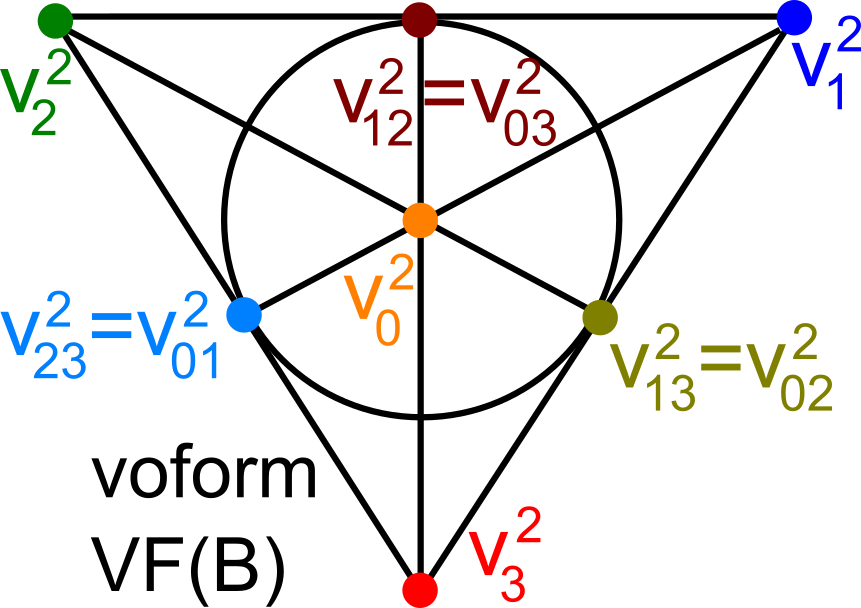}
\includegraphics[height=27mm]{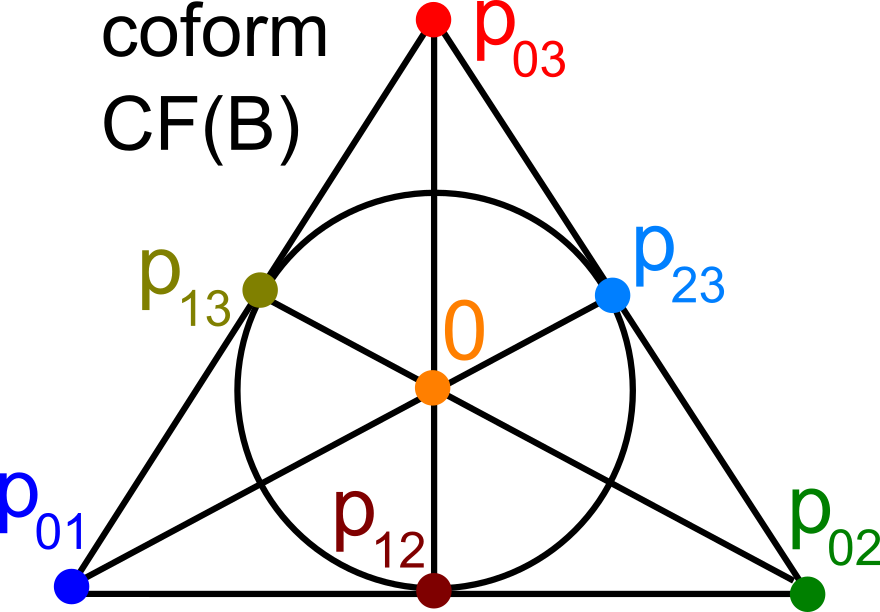}
\label{fig:forms3d}
\end{figure}

For a lattice $\La\subset\R^3$ with an obtuse superbase $B$, Definition~\ref{dfn:forms3d} introduces the voform $\VF(B)$ and the coform $\CF(B)$, which will be converted into root invariants later.
These forms are Fano planes marked by vonorms and conorms, respectively.
The \emph{Fano} projective plane of order 2 consists of seven non-zero classes (called \emph{nodes}) of the space $\La/2\La$, arranged in seven triples (called \emph{lines}).
If we mark these nodes by binary numbers $001$, $010$, $011$, $100$, $101$, $110$, $111$, the digit-wise sum of any two numbers in each line equals the third number modulo 2, see Fig.~\ref{fig:forms3d}.

\begin{dfn}[voform $\VF(B)$ and coform $\CF(B)$ of an obtuse superbase]
\label{dfn:forms3d} 
The \emph{voform} $\VF(B)$ of any obtuse superbase $B=(v_0,v_1,v_2,v_3)$ in $\R^3$ is the Fano plane in Fig.~\ref{fig:forms3d} with four nodes marked by $v_0^2,v_1^2,v_2^2,v_3^2$ and three nodes marked by $v_{12}^2, v_{23}^2,v_{13}^2$ so that $v_0^2$ is in the centre, $v_1^2$ is opposite to $v_{23}^2$, etc.
The \emph{coform} $\CF(B)$ is the dual Fano plane in Fig.~\ref{fig:forms3d} with three nodes marked by $p_{12},p_{23},p_{13}$ and three nodes marked by $p_{01},p_{02},p_{03}$, the centre is marked by $0$.
\bs
\end{dfn}

Much earlier than \cite{conway1992low}, Delone represented an obtuse superbase $B$ of $\vec a,\vec b, \vec c$, $\vec d=-\vec a-\vec b-\vec c$ by the skeleton of a tetrahedron with six (negative) scalar products on edges.
This \emph{Delone tetrahedron} is equivalent to the coform $\CF(B)$, which will be written in a matrix form in Definition~\ref{dfn:index-permutations}.
In 1975 \cite[chapter 10.4, p.~154]{delone1975bravais} claimed (without proof) a unique description of any lattice up to isometry by a 6-parameter Delone symbol satisfying sophisticated systems of equations and inequalities in 16 cases.
Theorem~\ref{thm:classification3d} will give a simpler and proved solution by root invariants in Definition~\ref{dfn:RI} based on only five Voronoi types.
\medskip

The \emph{zero} conorm $p_0=0$ at the centre of the coform $\CF(B)$ seems mysterious, because Conway and Sloane \cite{conway1992low} gave no formula for $p_0$, which also wrongly became non-zero in their Fig.~5.
This past mystery is explained by Lemma~\ref{lem:vonorms<->conorms}.

\begin{lem}[6 conorms $\lra$ 7 vonorms]
\label{lem:vonorms<->conorms}
For distinct indices $i,j\in\{0,1,2,3\}$, the conorm $p_{ij}$ in $\CF(B)$ of any superbase $B$ defines the dual line in the voform $\VF(B)$ through the nodes marked by $v_{ij}^2,v_k^2,v_l^2$ for $\{k,l\}=\{0,1,2,3\}-\{i,j\}$. 
Then
$$4p_{ij}=v_i^2+v_j^2+v_{ik}^2+v_{jk}^2-v_{ij}^2-v_k^2-v_l^2, \leqno{(\ref{lem:vonorms<->conorms}a)}$$
where the vonorms with negative signs are in the line of the voform $\VF(B)$ dual to $p_{ij}$.
The zero conorm $p_0=0$ in $\CF(B)$ can be computed by the similar formula
$$4p_{0}=v_0^2+v_1^2+v_2^2+v_3^2-v_{01}^2-v_{02}^2-v_{03}^2=0, \leqno{(\ref{lem:vonorms<->conorms}b)}$$
where the line dual to the zero conorm $p_0$ is the `circle' through $v_{01}^2,v_{02}^2,v_{03}^2$.
\bs
\end{lem}
\begin{proof}
Since all indices $i,j,k,l\in\{0,1,2,3\}$ are distinct, formula (\ref{lem:vonorms<->conorms}a) is symmetric in $k,l$ due to $v_{ik}^2+v_{jk}^2=v_{il}^2+v_{jl}^2$ following from $v_{ik}=v_i+v_k=-(v_j+v_l)=-v_{jl}$, $v_{jk}=v_j+v_k=-(v_i+v_l)=-v_{il}$.
To prove (\ref{lem:vonorms<->conorms}a), simplify its right hand side:
$$v_i^2+v_j^2+v_{ik}^2+v_{jk}^2-v_{ij}^2-v_k^2-v_l^2
=v_i^2+v_j^2+(v_i+v_k)^2+(v_j+v_k)^2-(v_i+v_j)^2-v_k^2-$$
$$-(-v_i-v_j-v_k)^2
=v_i^2+v_j^2+(v_i^2+2v_iv_k+v_k^2)+(v_j^2+2v_jv_k+v_k^2)-$$
$$-(v_i^2+2v_iv_j+v_j^2)-v_k^2-(v_i^2+v_j^2+v_k^2+2v_i v_j+2v_iv_k+2v_jv_k)=-4v_iv_j=4p_{ij}
.$$
(\ref{lem:vonorms<->conorms}b) follows from
$v_0^2+v_1^2+v_2^2+v_3^2=v_{01}^2+v_{02}^2+v_{03}^2$
in Definition~\ref{dfn:forms3d}.
\ws
\end{proof}
 
\begin{figure}
\label{fig:forms3d_permutations}
\caption{Actions of permutations $1\lra 2$ and $0\lra 1$ on voforms (top) and coforms.}
\includegraphics[width=\textwidth]{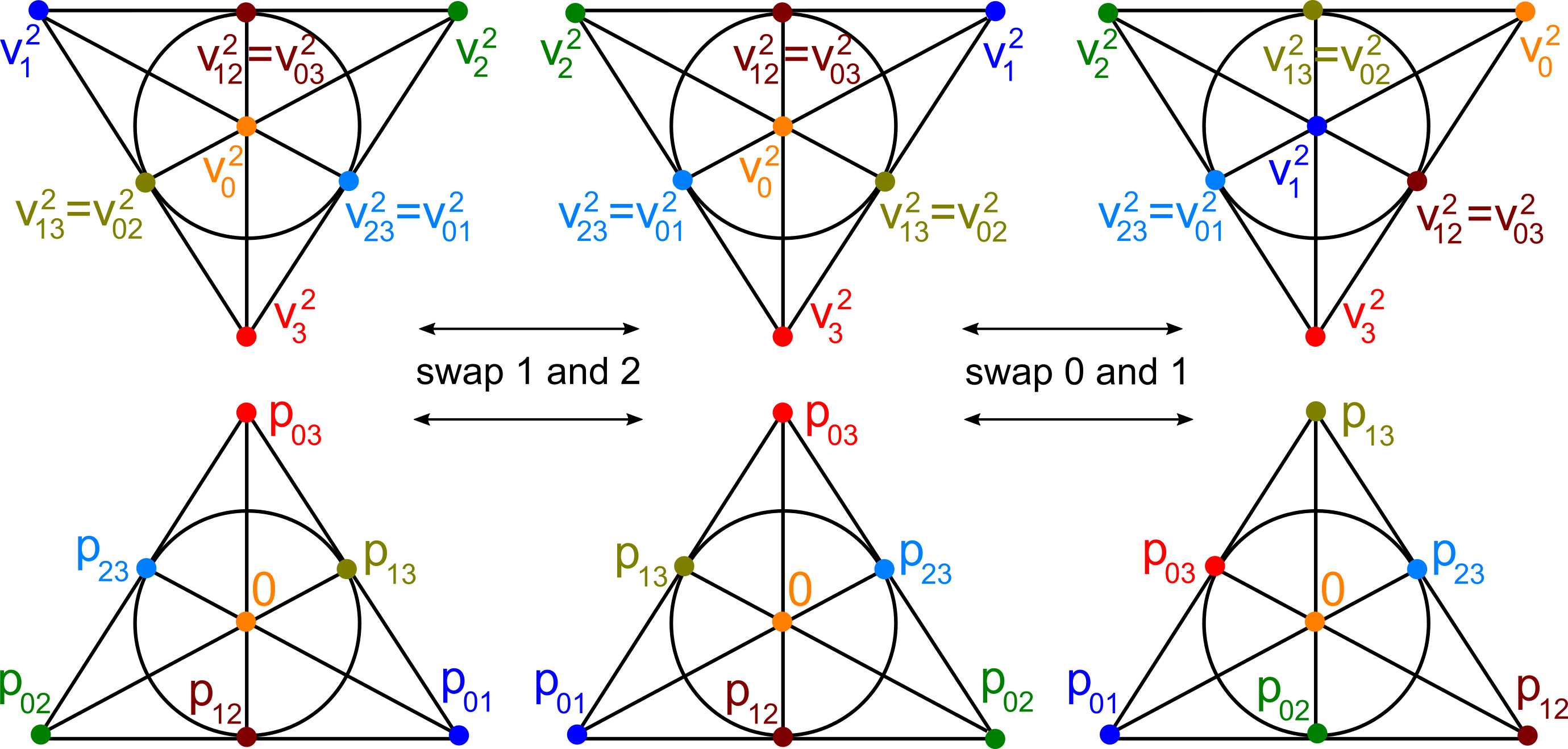}
\end{figure}

\begin{dfn}[index-permutations on vonorms and conorms]
\label{dfn:index-permutations}
For any \\ ordered obtuse superbase $B=\{v_0,v_1,v_2,v_3\}$, an \emph{index-permutation} is a permutation $\si\in S_4$ of indices $0,1,2,3$, which maps vonorms as follows: $v_i^2\mapsto v_{\si(i)}^2$, $v_{ij}^2\mapsto v_{\si(i)\si(j)}^2$, where $v_{ij}^2=v_{ji}^2$.
If we swap $v_1^2,v_2^2$, then we also swap only $v_{13}^2=v_{02}^2$ and $v_{23}^2=v_{01}^2$.
If we swap $v_0^2,v_1^2$, then we also swap only $v_{12}^2=v_{03}^2$ and $v_{02}^2=v_{13}^2$, see Fig.~\ref{fig:forms3d_permutations}.
Any index-permutation $\si\in S_4$ maps conorms by $p_{ij}\mapsto p_{\si(i)\si(j)}$, where $p_{ij}=p_{ji}$.
The group $S_4$ of all 24 index-permutations is generated by the three \emph{index-transpositions} $0\lra 1$, $1\lra 2$, $2\lra 3$.
\bs 
\end{dfn}

For any ordered superbase $\{v_0,v_1,v_2,v_3\}$, the voform can be written as the $2\times 3$ matrix $\VF(B)=\mat{v_{23}^2}{v_{13}^2}{v_{12}^2}{v_1^2}{v_2^2}{v_3^2}$, where $v_{23}^2=v_{01}^2$ is above $v_{1}^2$ and so on.
The 7th vonorm can be found as $v_0^2= v_{23}^2 + v_{13}^2 + v_{12}^2  - v_1^2 - v_2^2- v_3^2$ and is unnecessary to include into the matrix.
A coform can be written as 
$\CF(B)=\mat{p_{23}}{p_{13}}{p_{12}}{p_{01}}{p_{02}}{p_{03}}$.

\begin{lem}
\label{lem:index-permutations}
For any ordered obtuse superbase $B=\{v_0,v_1,v_2,v_3\}$, all 24 index-permutations act on the coform $\CF(B)$ as compositions of the transpositions: 
\smallskip

\noindent
\textbf{(a)}
$i\lra j$ for non-zero $i\neq j$ swaps the columns $i,j$ in $\CF(B)$, for example
$$\mat{\pac}{\pbc}{p_{12}}{\pdb}{\pda}{p_{03}}
\stackrel{1\lra 2}{\longleftrightarrow}
\mat{\pbc}{\pac}{\pab}{\pda}{\pdb}{\pdc}
\stackrel{0\lra 1}{\longleftrightarrow}
\mat{p_{23}}{\pdc}{\pdb}{p_{01}}{\pab}{\pac};
\leqno{(\ref{lem:index-permutations}a)}$$
\noindent
\textbf{(b)}
$0\lra i$ for $i\neq 0$ diagonally swaps pairs in the columns of indices $j\neq i,0$.
\medskip

\noindent
Any even permutation from $A_4$ acts as composition of the following permutations:
$$\mat{p_{23}}{\pdb}{\pdc}{p_{01}}{\pac}{\pab}
\stackrel{0\lra 1, 2\lra 3}{\longleftrightarrow}
\mat{\pbc}{\pac}{\pab}{\pda}{\pdb}{\pdc}
\stackrel{0\mapsto 1\mapsto 2\mapsto 0}{\longleftrightarrow}
\mat{\pdc}{\pbc}{\pdb}{\pab}{\pda}{\pac}.
\leqno{(\ref{lem:index-permutations}b)}$$
\end{lem}
\begin{proof}
By Definition~\ref{dfn:index-permutations} the action of any index-permutation $\si\in S_4$ on $\CF(B)$ follows by permuting the indices of conorms: $p_{ij}\mapsto p_{\si(i)\si(j)}$.
In all cases, any two conorms from one column of $\CF(B)$ remain in one column.
The composition of two transpositions such as $0\lra 1, 2\lra 3$ vertically swaps conorms in columns $2,3$.
The even permutation $1\mapsto 2\mapsto 3\mapsto 1$ cyclically permutes columns $1,2,3$.
Another even permutation $0\mapsto 1\mapsto 2\mapsto 0$ involving index 0 cyclically permutes the triples $(p_{23},p_{13},p_{03})$ of the coforms all including index 3 and the triple $(p_{01},p_{02},p_{12})$ of the coforms all excluding index 3.
\ws
\end{proof}

Lemma~\ref{lem:index-permutations} shows that coforms of six conorms are easier than voforms, which essentially require seven vonorms since $v_0^2$ appears after the transposition $0\lra 1$. 
$$\mat{
\color[rgb]{0,0.5,0.5}{v_{13}^2}}{
\color[rgb]{0,0.5,1}{v_{23}^2}}{
v_{12}^2}{
\color[rgb]{0,0.5,0}{v_2^2}}{
\color[rgb]{0,0,1}{v_1^2}}{
v_3^2}
\stackrel{1\lra 2}{\longleftrightarrow}\VF(B)=
\mat{
\color[rgb]{0,0.5,1}{v_{23}^2}}{
\color[rgb]{0,0.5,0.5}{v_{13}^2}}{
\color[rgb]{0.5,0,0}{v_{12}^2}}{
\color[rgb]{0,0,1}{v_1^2}}{
\color[rgb]{0,0.5,0}{v_2^2}}{
v_3^2}
\stackrel{0\lra 1}{\longleftrightarrow}
\mat{
v_{23}^2}{
\color[rgb]{0.5,0,0}{v_{12}^2}}{
\color[rgb]{0,0.5,0.5}{v_{13}^2}}{
\color[rgb]{1,0.5,0}{v_0^2}}{
v_2^2}{
v_3^2}$$

\begin{dfn}[odd-sum and even-sum vectors, digital sums]
\label{dfn:digital_sums}
For any basis $v_1,v_2,v_3$ in $\R^3$, write the partial sums $v_S$ from Lemma~\ref{lem:partial_sums} in coordinates:
\smallskip

\noindent
(\ref{dfn:digital_sums}o) 
$v_1=(1,0,0)$,
$v_2=(0,1,0)$, 
$v_3=(0,0,1)$, 
$v_0
=(-1,-1,-1)$;
\smallskip

\noindent
(\ref{dfn:digital_sums}e) 
$v_{12}=v_1+v_2=(1,1,0)$,
$v_{23}=v_2+v_3=(0,1,1)$,
$v_{13}=v_1+v_3=(1,0,1)$.
\smallskip

\noindent
The four vectors (and their opposites) from (\ref{dfn:digital_sums}o) are called \emph{odd-sum} vectors, because the sum of their coordinates is odd.
The three vectors (and their opposites) from (\ref{dfn:digital_sums}e) are called \emph{even-sum} vectors.
For any vector $v=(x_1,x_2,x_3)$ with coordinates $x_1,x_2,x_3\in\Z$, its \emph{digital} image is $[v]=100x_1+10x_2+x_3$.
\bs
\end{dfn}

\begin{lem}[digital sums sufficiency]
\label{lem:digital_sums}
For any basis $v_1,v_2,v_3$, let $u,v$ be sums of at most four
vectors $v_S$ in Lemma~\ref{lem:partial_sums}.
Then $u=v$ if and only if $[u]=[v]$.
\bt
\end{lem}
\begin{proof}
Any partial sum or its opposite from (\ref{dfn:digital_sums}o,e) has all coordinates in the range $[-1,1]$.
Both $u,v$ have coordinates in the range $[-4,4]$.
The equality between $[v]=100x_1+10x_2+x_3$ and $[u]=100y_1+10y_2+y_3$ is equivalent to $100(x_1-x_2)+10(y_1-y_2)+(z_1-z_2)=0$.
Since each integer difference in the brackets is within $[-8,8]$, the last equality can hold only if all differences vanish, so $u=v$.
\ws
\end{proof}

\section{An explicit description all obtuse superbases of 3D lattices}
\label{sec:superbases}

In this section Lemmas~\ref{lem:V1superbases}-\ref{lem:V5superbases} describe all possible obtuse superbases of any lattice $\La\subset\R^3$.
Any obtuse superbase $\{v_0,v_1,v_2,v_3\}$ has its dual $\{-v_0,-v_1,-v_2,-v_3\}$ related by the central symmetry with respect to $0$.
Lemmas~\ref{lem:V1superbases}-\ref{lem:V5superbases} describe all obtuse superbases and their coforms separately for each Voronoi type.
\medskip

Even in the generic case, \cite[chapter 7.5, p.~130]{delone1975bravais} missed the following step and went straight to Delone paramaters of a single obtuse superbase.
Lemma~\ref{lem:V1superbases} proves that any lattice $\La$ of a Voronoi type $V_1$ has only one pair centrally symmetric obtuse superbases.
There will be more non-isometric obtuse superbases for higher symmetry types in Lemmas~\ref{lem:V2superbases}-\ref{lem:V5superbases}.


\begin{lem}[obtuse superbases for Voronoi type $V_1$]
\label{lem:V1superbases} 
Let a lattice $\La\subset\R^3$ have Voronoi type $V_1$, so the Voronoi domain $V(\La)$ is a truncated octahedron.
\smallskip

\noindent
\textbf{(a)}
$\La$ has two obtuse superbases related by the central symmetry $v\lra -v$;
\smallskip

\noindent
\textbf{(b)}
coforms of all obtuse superbases of $\La$ are related by 24 index-permutations.
\bt 
\end{lem}
\begin{proof}
\textbf{(a)}
Let $\{v_0,v_1,v_2,v_3\}$ be any obtuse superbase of $\La$, which exists by Theorem~\ref{thm:reduction}.
Since the Voronoi domain $V(\La)$ is a truncated octahedron with seven pairs of parallel opposite faces.
The lattice $\La$ has seven pairs of strict Voronoi vectors orthogonal to these faces.
By Lemma~\ref{lem:partial_sums} all these seven pairs of Voronoi vectors should coincide with the partial sums and their opposites $\pm v_S$ from (\ref{dfn:digital_sums}o,e).
\medskip

The Voronoi domain $V(\La)$ has four pairs of opposite hexagonal faces obtained by cutting corners in four pairs of opposite triangular faces of an octahedron.
The normal vectors of these hexagons are Voronoi odd-sum vectors $\pm v_i$, $i=0,1,2,3$.
The Voronoi even-sum vectors $v_{ij}=v_i+v_j$ 
are normal to the three pairs of opposite parallelogram faces obtained by cutting three pairs of opposite vertices.
\medskip

The seven pairs of Voronoi vectors have these digital sums from Definition~\ref{dfn:digital_sums}:

\noindent
(\ref{lem:V1superbases}o) 
Voronoi odd-sums 
$[\pm v_1]=\pm 100$,
$[\pm v_2]=\pm 10$,
$[\pm v_3]=\pm 1$,
$[\pm v_0]=\mp 111$.
\smallskip

\noindent
(\ref{lem:V1superbases}e) 
Voronoi even-sum vectors 
$[\pm v_{12}]=\pm 110$,
$[\pm v_{23}]=\pm 11$, 
$[\pm v_{13}]=\pm 101$.
\medskip

If an obtuse superbase $\{u_0,u_1,u_2,u_3\}$ consists of four odd-sum vectors, by Lemma~\ref{lem:digital_sums} the condition $u_0+u_1+u_2+u_3=0$ is equivalent to $[u_0]+[u_1]+[u_2]+[u_3]=0$ for some digital sums from (\ref{lem:V1superbases}o).
The only possibility $100+10+1+(-111)=0$ up to a sign gives the known obtuse superbases $\pm\{v_0,v_1,v_2,v_3\}$.
If an obtuse superbase has one even-sum vector $u_0$, then it should have one more, say $u_1$, otherwise an odd sum $[u_1]+[u_2]+[u_3]$ cannot become 0 after adding an even integer $[u_0]$.  
For any choice of $u_0\neq\pm u_1$ from (\ref{lem:V1superbases}e), by Lemma~\ref{lem:partial_sums} the sum $u_0+u_1$ should be another even-sum vector from (\ref{lem:V1superbases}e).
But there is no choice of signs such that $\pm 110\pm 11\pm 101=0$. 
\medskip

\noindent
\textbf{(b)}
By part (a) all obtuse superbases $B$ of $\La$ differ either by re-ordering vectors of $B$ or by the central symmetry with respect to the origin of $\R^3$, which keeps the coform $\CF(B)$ invariant.
Lemma~\ref{lem:index-permutations}(c) implies that coforms $\CF(B)$ of all obtuse superbases $B$ of $\La$ are related by 24 index-permutations from Definition~\ref{dfn:index-permutations}.
\ws
\end{proof}

\begin{lem}[obtuse superbases for Voronoi type $V_2$]
\label{lem:V2superbases} 
Let a lattice $\La\subset\R^3$ have Voronoi type $V_2$, so the Voronoi domain is a hexa-rhombic dodecahedron.
\smallskip

\noindent
\textbf{(a)}
$\La$ has an obtuse superbase $\{v_0,v_1,v_2,v_3\}$ with one pair of orthogonal vectors, say $v_2\cdot v_3=0$.
Then any obtuse superbase $B$ of $\La$ is isometric to one of
$$\text{obtuse superbases }B_1=\{v_0,v_1,v_2,v_3\} \text{ and } 
B_2=\{v_0+v_3, v_1+v_3, v_2, -v_3\}.
\leqno{(\ref{lem:V2superbases})}$$

\noindent
\textbf{(b)} 
Let an obtuse superbase $B=\{v_0,v_1,v_2,v_3\}$ with $v_1\cdot v_2=0$ have a coform $\CF(B_1)=\mat{0}{\pac}{p_{12}}{p_{01}}{p_{02}}{\pdc}$ with $p_{23}=0$.
Then another obtuse superbase $B_2=\{v_0+v_3, v_1+v_3, v_2, -v_3\}$ has the coform $\CF(B_2)=\mat{0}{\pdc}{p_{12}}{p_{01}}{p_{02}}{\pac}$.
\medskip

\noindent
\textbf{(c)} 
Any obtuse superbase $B$ of $\La$ has exactly one zero conorm.
The 24 index-permutations from Definition~\ref{dfn:index-permutations} allow us to write $\CF(B)=\mat{0}{\pac}{\pab}{p_{01}}{\pdb}{\pdc}$.
The above forms with $p_{23}=0$ over all obtuse superbases of $\La$ are related by the symmetry group $D_4$ (of a square) acting on the $2\times 2$ submatrix $\matfour{\pac}{\pab}{\pdb}{\pdc}$. 
\bt 
\end{lem}
\begin{proof}
\textbf{(a)}
In comparison with the most generic Voronoi domain in Lemma~\ref{lem:V1superbases}, a hexa-rhombic dodecahedron $V(\La)$ has six pairs of faces: two pairs of hexagons and four pairs parallelograms with one pair degenerated from one pair of hexagons and one pair of parallelograms disappeared.
This degeneracy appears when exactly two of four superbase vectors become orthogonal, say $v_2\cdot v_3=0$.
\medskip

In addition to the seven pairs of Voronoi vectors $\pm v_S$ from Lemma~\ref{lem:V1superbases}, we have exactly one extra pair of the non-strict Voronoi vectors $\pm(v_3-v_2)$ whose length equals $|v_3+v_2|$ due to $v_2\cdot v_3=0$.
Now we have the extra choice of the Voronoi even-sum vector $v_3-v_2=(0,-1,1)$ and its opposite.
Here are the digital images of all $4\times 2+4\times 2$ Voronoi vectors.
\smallskip

\noindent
(\ref{lem:V2superbases}o) 
Voronoi odd-sums
$[\pm v_1]=\pm 100$,
$[\pm v_2]=\pm 10$,
$[\pm v_3]=\pm 1$,
$[\pm v_0]=\mp 111$.
\smallskip

\noindent
(\ref{lem:V2superbases}e) 
even : 
$[\pm v_{12}]=\pm 110$,
$[\pm v_{23}]=\pm 11$, 
$[\pm v_{13}]=\pm 101$,
$[\pm(v_3-v_2)]=\mp 9$.
\medskip

If an obtuse superbase has only four odd-sum vectors from (\ref{lem:V2superbases}e), we can get only $\pm B_1=\pm\{v_0,v_1,v_2,v_3\}$ as in Lemma~\ref{lem:V1superbases}.
Choosing an even-sum vector $u_0$ from (\ref{lem:V2superbases}e), we should include at least one more even-sum vector, say $u_1$.
The negative partial sum $-u_0-u_1$ by Lemma~\ref{lem:partial_sums} should be among other Voronoi even-sum vectors in (\ref{lem:V2superbases}e) so that $[u_0]+[u_1]+[-u_0-u_1]=0$.
Without a vector from the new pair $\pm(v_3-v_2)$, no choice of signs gives $0=\pm 110\pm 11\pm 101$.
\medskip

The only possible identity $[u_0]+[u_1]+[-u_0-u_1]=0$ with a new digital sum $9$ from (\ref{lem:V2superbases}e) is $110-101-9=0$ up to a permutation and an overall sign.
Hence we can get another obtuse superbase (potentially not isometric to $B_1$) only by choosing $u_0=-v_{12}=-v_1-v_2=v_0+v_3$ with $[u_0]=-110$ and $u_1=v_{13}=v_1+v_3$ with $[u_1]=101$ so that $u_0+u_1=v_3-v_2$ with $[u_0]=-9$ (up to a sign and re-ordering). 
\medskip

Other superbase vectors $u_2,u_3$ should have the digital sum $[u_2]+[u_3]=-[u_0]-[u_1]=110-101=9$.
The remaining digital sums from (\ref{lem:V2superbases}o) and~(\ref{lem:V2superbases}e) give only one splitting $9=10-1$, so $u_2=v_2$, $u_3=-v_3$.  
We got the second obtuse superbase from (\ref{lem:V2superbases}): $B_2=\{u_0,u_1,u_2,u_3\}=\{v_0+v_3, v_1+v_3, v_2, -v_3\}$.
\medskip

Using the transposition $2\lra 3$ of indices and respecting $v_2\cdot v_3=0$, we get the obtuse superbase 
$\{v_0+v_2,v_1+v_2,v_3,-v_2\}=\{-v_1-v_3,-v_0-v_3,v_3,-v_2\}$.
After re-ordering, the last superbase becomes opposite (isometric via $v\mapsto -v$) to $B_2$.
Under the transposition $0\lra 1$, vectors are only permuted in both $B_1,B_2$.
\medskip 

\noindent
\textbf{(b)}
The conorms $q_{ij}$ of the superbase $B_2=\{v_0+v_3, v_1+v_3, v_2, -v_3\}$ are expressed via the conorms $p_{ij}$ of $B_1=\{v_0,v_1,v_2,v_3\}$ with $p_{23}=-v_2\cdot v_3=0$ as follows:

$q_{23}=-v_2\cdot(-v_3)=v_2\cdot v_3=0$,

$q_{13}=-(v_1+v_3)\cdot(-v_3)=(v_0+v_2)\cdot(-v_3)=-v_0\cdot v_3=p_{03}$,

$q_{12}=-(v_1+v_3)\cdot v_2=-v_1\cdot v_2=p_{12}$, 

$q_{01}=-(v_0+v_3)\cdot(v_1+v_3)=(v_1+v_2)\cdot(v_1+v_3)=v_1(v_1+v_2+v_3)=-v_1\cdot v_0=p_{01}$,

$q_{02}=-(v_0+v_3)\cdot v_2=-v_0\cdot v_2=p_{02}$,

$q_{03}=-(v_0+v_3)\cdot(-v_3)=(v_1+v_2)\cdot(-v_3)
=p_{13}$, so $\CF(B_2)=\mat{0}{\pdc}{p_{12}}{p_{01}}{p_{02}}{\pac}$.
\medskip

\noindent
\textbf{(c)}
By parts \textbf{(a,b)} and Lemma~\ref{lem:index-permutations}, the coform of any superbase of $\La$ up to 24 index-permutations is either 
$\CF(B_1)=\mat{0}{\pac}{p_{12}}{p_{01}}{p_{02}}{\pdc}$ or 
$\CF(B_2)=\mat{0}{\pdc}{p_{12}}{p_{01}}{p_{02}}{\pac}$ 
related by the transposition $\pac\lra\pdc$.
Keeping the first column fixed, the index-permutations induced by $2\lra 3$ and $0\lra 1$ act on $\matfour{\pac}{\pab}{\pdb}{\pdc}$ by swapping the columns and by swapping diagonally opposite elements, see Definition~\ref{dfn:index-permutations}.
\medskip

\noindent
Use $\pmb{\CF(B_1)}$: 
$\matfour{\pab}{\pac}{\pdc}{\pdb}\stackrel{2\lra 3}{\mapsfrom}
\pmb{\matfour{\pac}{\pab}{\pdb}{\pdc}}\stackrel{0\lra 1}{\mapsto}
\matfour{\pdc}{\pdb}{\pab}{\pac}\stackrel{2\lra 3}{\mapsto}
\matfour{\pdb}{\pdc}{\pac}{\pab}$.

\noindent
Use $\pmb{\CF(B_2)}$: 
$\matfour{\pab}{\pdc}{\pac}{\pdb}\stackrel{2\lra 3}{\mapsfrom}
\pmb{\matfour{\pdc}{\pab}{\pdb}{\pac}}\stackrel{0\lra 1}{\mapsto}
\matfour{\pac}{\pdb}{\pab}{\pdc}\stackrel{2\lra 3}{\mapsto}
\matfour{\pdb}{\pac}{\pdc}{\pab}$.
The eight arrangements above are realised by 
the symmetry group $D_4$ of a square. 
\end{proof}

\begin{lem}[obtuse superbases for Voronoi type $V_3$]
\label{lem:V3superbases} 
Let a lattice $\La\subset\R^3$ have Voronoi type $V_3$, so the Voronoi domain $V(\La)$ is a rhombic dodecahedron.
\medskip

\noindent
\textbf{(a)}
In this case the lattice $\La$ 
has an obtuse superbase $B_1=\{v_0,v_1,v_2,v_3\}$ with two different pairs of orthogonal vectors, say $v_0\cdot v_1=0=v_2\cdot v_3$.
Then any obtuse superbase of $\La$ is isometric to one of the following obtuse superbases:
$$B_1,\quad
B_2=\{v_0+v_3, v_1+v_3, v_2, -v_3\},\quad
B_3=\{v_0,-v_1,v_2+v_1,v_3+v_1\}.
\leqno{(\ref{lem:V3superbases})}$$

\noindent
\textbf{(b)}
Any obtuse superbase $B$ of $\La$ has exactly two zero conorms in one column.
The 24 index-permutations from Definition~\ref{dfn:index-permutations} allow us to write $\CF(B)=\mat{0}{\pac}{\pab}{0}{\pdb}{\pdc}$.
In the above form with $p_{23}=0=p_{01}$ for all obtuse superbases of $\La$, the four non-zero conorms can be freely permuted by the symmetry group $S_4$. 
\bt
\end{lem}
\begin{proof}
\textbf{(a)}
In comparison with Lemma~\ref{lem:V2superbases}, a rhombic dodecahedron has five pairs of parallelograms degenerated from a 12-face hexa-rhombic dodecahedron due to another pair of orthogonal vectors, say $v_0\cdot v_1=0$ in addition to $v_2\cdot v_3=0$.
\medskip

This degeneracy adds the 9th pair of non-strict Voronoi vectors $\pm(v_0-v_1)$ whose lengths equals $|v_0+v_1|$ since $v_0\cdot v_1=0$.
The first two superbases $B_1,B_2$ in (\ref{lem:V3superbases}) came from Lemma~\ref{lem:V2superbases}.
The double transposition of indices $0\lra 2$, $1\lra 3$ respects $v_0\cdot v_1=0=v_2\cdot v_3$ and maps $B_2=\{v_0+v_3, v_1+v_3, v_2, -v_3\}$ as follows:
$$B_2\stackrel{0\lra 2, 1\lra 3}{\mapsto}\{v_2+v_1,v_1+v_3,v_0,-v_1\}=
\{v_0,-v_1,v_2+v_1,v_3+v_1\}=B_3.$$

We will show that any other obtuse superbase is isometric to one of $B_1,B_2,B_3$.
We have four pairs of odd-sum vectors and five pairs of even-sum vectors below:
\smallskip

\noindent
(\ref{lem:V3superbases}o) 
Voronoi odd-sums 
$[\pm v_1]=\pm 100$,
$[\pm v_2]=\pm 10$,
$[\pm v_3]=\pm 1$,
$[\pm v_0]=\mp 111$;
\smallskip

\noindent
(\ref{lem:V3superbases}e) 
Voronoi even-sum vectors 
$[\pm v_{12}]=\pm 110$,
$[\pm v_{23}]=\pm 11$, 
$[\pm v_{13}]=\pm 101$,
and
$[\pm(v_3-v_2)]=\mp 9$,
$[\pm(v_0-v_1)]=\mp 211$.
\medskip

Since the condition $v_0\cdot v_1=0$ wasn't used in Lemma~\ref{lem:V2superbases}, it suffices to consider only superbases whose partial sums have a vector from the new pair $\pm(v_0-v_1)$.
\medskip

Looking for even-sum vectors $u_0,u_1$ and $-u_0-u_1$ from (\ref{lem:V2superbases}e), the only possible identity $[u_0]+[u_1]+[-u_0-u_1]=0$ with a new digital sum $211$ from (\ref{lem:V2superbases}e) is $211-110-101=0$ up to a permutation and sign.
We can get another obtuse superbase not isometric to $B_1,B_2$ from Lemma~\ref{lem:V2superbases} only by choosing $u_0=-v_{12}=-v_1-v_2=v_0+v_3$, $[u_0]=-110$ and 
$u_1=-v_{13}=-v_1-v_3=v_0+v_2$, $[u_1]=-101$ so that $u_0+u_1=-v_1-v_2+v_0+v_2=v_0-v_1$, $[u_0+u_1]=-211$ (up to a sign and re-ordering). 
Other superbase vectors $u_2,u_3$ should have the digital sum $[u_2]+[u_3]=-[u_0]-[u_1]=211$.
The remaining digital sums from (\ref{lem:V2superbases}o) and~(\ref{lem:V2superbases}e) give only $211=111+100$, so $u_2=-v_0$, $u_1=v_1$.  
This superbase $\{-v_1-v_2,-v_3-v_3,-v_0,v_1\}$ is opposite to $B_3=\{v_0,-v_1,v_2+v_1,v_3+v_1\}$ up to re-ordering.
\medskip

\noindent
\textbf{(b)}
If $B_1=\{v_0,v_1,v_2,v_3\}$ with $v_0\cdot v_1=0=v_2\cdot v_3$ has a coform 
$\CF(B_1)=\mat{0}{\pac}{\pab}{0}{\pdb}{\pdc}$
with $p_{12}=0=p_{03}$, 
the superbase $B_2=\{v_0+v_3, v_1+v_3, v_2, -v_3\}$ has 
$\CF(B_2)=\mat{0}{\pdc}{p_{12}}{0}{p_{02}}{\pac}$ by
Lemma~\ref{lem:V2superbases}(b) with the extra restriction $p_{01}=0$.
\medskip

The obtuse superbase $B_3=\{v_0,-v_1,v_2+v_1,v_3+v_1\}$ potentially non-isometric to $B_1,B_2$ has the conorms $q_{ij}$ expressed via the conorms $p_{ij}$ of $B_1$ as follows:

$q_{23}=-(v_2+v_1)\cdot(v_3+v_1)=-v_1\cdot (v_1+v_2+v_3)=v_1\cdot v_0=0$,

$q_{13}=v_1\cdot(v_3+v_1)=-v_1\cdot (v_0+v_2)=-v_1\cdot v_2=p_{12}$, 

$q_{12}=v_1\cdot(v_2+v_1)=-v_1\cdot (v_0+v_3)=-v_1\cdot v_3=p_{13}$,

$q_{01}=-v_0\cdot(-v_1)=0$,

$q_{02}=-v_0\cdot(v_2+v_1)=-v_0\cdot v_2=p_{02}$,

$q_{03}=-v_0\cdot(v_3+v_1)=-v_0\cdot v_3=p_{03}$.
Hence $\CF(B_3)=\mat{0}{\pab}{\pac}{0}{p_{02}}{p_{03}}$.
\medskip

\noindent
By part \textbf{(a)} and Lemma~\ref{lem:index-permutations}, the coform of any superbase of $\La$ up to 24 index-permutations is one of 
the above coforms $\CF(B_1),\CF(B_2),\CF(B_3)$, which are
related by the transpositions $\pac\lra\pdc$ and $\pab\lra\pac$.
Keeping the first column fixed in a coform, the index-permutations induced by $2\lra 3$ and $0\lra 1$ act on $\matfour{\pac}{\pab}{\pdb}{\pdc}$ by swapping the columns and by swapping diagonally opposite elements.
All permutations above generate the full group $S_4$ permuting all non-zero conorms in the submatrix $\matfour{\pbc}{\pda}{\pac}{\pdb}$ of $\CF(B)$.
\ws
\end{proof}

\begin{lem}[obtuse superbases for Voronoi type $V_4$]
\label{lem:V4superbases} 
Let a lattice $\La\subset\R^3$ have Voronoi type $V_4$, so the Voronoi domain $V(\La)$ is a hexagonal prism.
\medskip

\noindent
\textbf{(a)}
Then $\La$ has an obtuse superbase
$B_1=\{v_0,v_1,v_2,v_3\}$ with one vector (say) $v_3$ orthogonal to two others $v_1,v_2$.
Any obtuse superbase of $\La$ is isometric to one of
$$
B_1,\quad
B_2=\{v_0+v_3, v_1+v_3, v_2, -v_3\},\quad
B_4=\{v_0+v_3, v_2+v_3, v_1, -v_3\}.
\leqno{(\ref{lem:V4superbases})}$$

\noindent
\textbf{(b)} 
Any obtuse superbase $B$ of $\La$ has exactly two zero conorms in one column.
The 24 index-permutations from Definition~\ref{dfn:index-permutations} allow us to write $\CF(B)=\mat{0}{0}{\pab}{\pda}{\pdb}{\pdc}$.
In the above form with $p_{23}=0=p_{13}$ for all obtuse superbases of $\La$, the conorms $\pab,\pda,\pdb$ can be freely permuted by the symmetry group $S_3$. 
\bt
\end{lem}
\begin{proof}
A hexagonal prism $V(\La)$ has four pairs of opposite parallel faces: one pair of hexagons and three pairs of rectangles.
So $V(\La)$ can be considered as a degenerate case of a hexa-rhombic dodecahedron from Lemma~\ref{lem:V2superbases}, not a rhombic dodecahedron.
This degeneracy happens due to one vector (say) $v_3$ orthogonal to other two superbase vectors $v_1,v_2$.
The first two superbases in (\ref{lem:V4superbases}) are inherited from  Lemma~\ref{lem:V2superbases}.
The last obtuse superbase in (\ref{lem:V4superbases}) is obtained from the second one by the transposition $1\lra 2$ of indices, which respects the new orthogonality conditions $v_1\cdot v_3=0=v_2\cdot v_3$. 
\medskip

The 2D lattice $\La_2$ has two pairs of obtuse superbases $\pm\{v_1,v_2,-v_1-v_2\}$.
We can choose any two of the three vectors $v_1,v_2,-v_1-v_2$ and complement this pair $u_2,u_3$ by $u_1=-v_3$ and $u_0=-u_1-u_2-u_3$.
The resulting three superbases are isometric to $B_1,B_2,B_4$ in (\ref{lem:V4superbases}).
For example, the superbase with $u_1=v_1$, $u_2=v_2$, $u_3=-v_3$, $u_0=-v_1-v_2+v_3$ is isometric to $B_1$
by $v_1\mapsto v_1$, $v_2\mapsto v_2$, $v_3\mapsto -v_3$.
\medskip

We check that any obtuse superbases of $\La$ is isometric to one of $B_1,B_2,B_4$.
In addition to the eight pairs of Voronoi vectors $\pm v_0,\pm v_1,\pm v_2,\pm v_3,\pm v_{12},\pm v_{23},\pm v_{13}$ and $\pm(v_3-v_2)$ in Lemma~\ref{lem:V2superbases}, we have two more pairs of non-strict Voronoi vectors $\pm(v_3-v_1)$ and $\pm(v_1+v_2-v_3)$ whose lengths are equal to $|v_3+v_1|$ and $|v_0|=|v_1+v_2+v_3|$, respectively, due to $v_1\cdot v_3=0=v_2\cdot v_3$.
We have 5+5 pairs: 
\smallskip

\noindent
(\ref{lem:V3superbases}o) 
Voronoi odd-sum vectors: 
$[\pm v_1]=\pm 100$,
$[\pm v_2]=\pm 10$,
$[\pm v_3]=\pm 1$, \\ and
$[\pm v_0]=\mp 111$,
$[\pm(v_1+v_2-v_3)]=\pm 109$;
\smallskip

\noindent
(\ref{lem:V3superbases}e) 
Voronoi even-sum vectors: 
$[\pm v_{12}]=\pm 110$,
$[\pm v_{23}]=\pm 11$, 
$[\pm v_{13}]=\pm 101$,
and
$[\pm(v_3-v_2)]=\mp 9$,
$[\pm(v_3-v_1)]=\mp 99$.
\medskip

Since Lemma~\ref{lem:V2superbases} didn't use the condition $v_1\cdot v_3=0$, it suffices to check only superbases whose partial sums have a vector from $\pm(v_3-v_1)$, $\pm(v_1+v_2-v_3)$.
\medskip

\textbf{Case of four odd-sum vectors}.
Trying to find four digital sums from (\ref{lem:V4superbases}e) to fit $[u_0]+[u_1]+[u_2]+[u_3]=0$, we conclude that one of $\pm 100$ and one of $\pm 10$ should be used, because three other pairs have odd digital sums.
Choosing one positive sign of (say) $100$, we have two sums $100\pm 10$.
The sum $90$ cannot be split as a sum of two numbers from $\{\pm 1,\pm 111,\pm 109\}$.
The only splitting $110=111-1$ of another sum misses $\pm 109$ and leads to the first superbase $B_1=\{v_0,v_1,v_2,v_3\}$.
\medskip

\textbf{Case of at least two even-sum vectors}.
Looking for even-sum vectors $u_0,u_1$ and $-u_0-u_1$ from (\ref{lem:V4superbases}e), the only possible identity $[u_0]+[u_1]+[-u_0-u_1]=0$ with a new digital sum $99$ 
is $110-11-99=0$ up to a sign and re-ordering.
\medskip

We can get another obtuse superbase not isometric to $B_1,B_2$ from Lemma~\ref{lem:V2superbases} only by choosing $u_0=-v_1-v_2=v_0+v_3$, $[u_0]=-110$ and 
$u_1=v_{23}=v_2+v_3$, $[u_1]=11$ so that $u_0+u_1=v_3-v_1$, $[u_0+u_1]=-99$ up to a sign and re-ordering. 
\medskip

Other vectors $u_2,u_3$ should have the digital sum $[u_2]+[u_3]=-[u_0]-[u_1]=99$.
The remaining digital sums from (\ref{lem:V4superbases}o) and~(\ref{lem:V4superbases}e) can give only one new splitting $99=100-1$, because we have already used $99=110-11$ above.
Choosing $u_2=v_1$ and $u_3=-v_3$ (up to a swap), we get $B_4=\{v_0+v_3,v_2+v_3,v_1,-v_3\}$ from (\ref{lem:V4superbases}). 
\medskip

\noindent
\textbf{(b)}
If $B_1=\{v_0,v_1,v_2,v_3\}$ with $v_1\cdot v_3=0=v_2\cdot v_3$ has a coform 
$\CF(B_1)=\mat{0}{0}{\pab}{p_{01}}{\pdb}{\pdc}$
with $p_{23}=0=p_{13}$, 
the superbase $B_2=\{v_0+v_3, v_1+v_3, v_2, -v_3\}$ has 
$\CF(B_2)=\mat{0}{\pdc}{\pab}{p_{01}}{\pdb}{0}$ by
Lemma~\ref{lem:V2superbases}(b) with the extra restriction $p_{13}=0$.
The index-permutation induced by $0\lra 1$ from Definition~\ref{dfn:index-permutations} transforms the above coform into $\CF(B_2)=\mat{0}{0}{\pdb}{p_{01}}{\pab}{\pdc}$ with $p_{23}=0=p_{13}$ as in $\CF(B_1)$.
\medskip

The obtuse superbase $B_4=\{v_0+v_3,v_2+v_3,v_1,-v_3\}$ potentially non-isometric to $B_1,B_2$ has the conorms $q_{ij}$ expressed via the conorms $p_{ij}$ of $B_1$ as follows:

$q_{23}=-v_1\cdot(-v_3)=0$,

$q_{13}=-(v_2+v_3)\cdot(-v_3)=(v_0+v_1)\cdot(-v_3)=-v_0\cdot v_3=p_{03}$,

$q_{12}=-(v_2+v_3)\cdot v_1=-v_2\cdot v_1=p_{12}$,

$q_{01}=-(v_0+v_3)\cdot(v_2+v_3)=(v_0+v_3)\cdot(v_0+v_1)=v_0(v_0+v_1+v_3)=-v_0\cdot v_2=p_{02}$,

$q_{02}=-(v_0+v_3)\cdot v_1=-v_0\cdot v_1=p_{01}$,

$q_{03}=-(v_0+v_3)\cdot(-v_3)=(v_1+v_2)\cdot(-v_3)=0$.

The index-permutation induced by $0\lra 1$ from Definition~\ref{dfn:index-permutations} transforms the resulting coform $\mat{0}{\pdc}{\pab}{\pdb}{\pda}{0}$ into $\CF(B_4)=\mat{0}{0}{\pda}{\pdb}{\pab}{\pdc}$.
Comparing $\CF(B_1)$, $\CF(B_2)$, $\CF(B_4)$, which all have $p_{23}=0=p_{13}$, we notice that $\pdc$ remains at the same place.
Actually, $p_{03}=-v_0\cdot v_3=(v_1+v_2+v_3)\cdot v_3=v_3^2$ is invariant as the squared length of the vector $v_3$ orthogonal to $v_2,v_3$. 
The other three conorms $\pda,\pdb,\pab$ have three arrangements in $\CF(B_1)$, $\CF(B_2)$, $\CF(B_4)$.
If we swap the first two columns by the index-permutation $1\lra 2$, we get all six arrangements.
Hence $\pab,\pda,\pdb$ can be freely permuted by the group $S_3$. 
\ws
\end{proof}

As an alternative to Lemma~\ref{lem:V4superbases} for any lattice $\La\subset\R^3$ of Voronoi type $V_4$, the index-permutation induced by $1\lra 2, 0\lra 3$ from Definition~\ref{dfn:index-permutations} re-writes $\CF(\La)$ as $\mat{\pab}{\pda}{\pdb}{0}{0}{\pdc}$ with fixed $\pdc$ and freely permutable conorms in the top row.

\begin{lem}[obtuse superbases for Voronoi type $V_5$]
\label{lem:V5superbases} 
Let a lattice $\La\subset\R^3$ have Voronoi type $V_4$, so the Voronoi domain $V(\La)$ is a cuboid.
Then any obtuse superbase of $\La$ belongs to one of the four isometry classes of obtuse superbases. 
\smallskip

\noindent
(\ref{lem:V5superbases}o) 
One class of 8 \emph{odd superbases} : 
$\{v_0,\pm v_1,\pm v_2,\pm v_3\}$ for any choice of signs.
Any coform can be written as $\mat{0}{0}{0}{|v_1|^2}{|v_2|^2}{|v_3|^2}$ up to 24 index-permutations.
\smallskip

\noindent
(\ref{lem:V5superbases}e) 
Three classes each consisting of 8 \emph{even superbases}
$\{v_i,v_j,v_k-v_i,-v_k-v_j\}$ for pairwisely orthogonal basis vectors $v_i,v_j,v_k\in\{\pm v_1,\pm v_2,\pm v_3\}$ with distinct $i,j,k\in\{1,2,3\}$.
Any coforms can be written as $\mat{0}{0}{|v_i|^2}{0}{|v_k|^2}{|v_j|^2}$ up to 24 index-permutations, 
where $k=1,2,3$ determines a class, but $i,j$ can be swapped.
\bt
\end{lem}
\begin{proof}
A cuboid $V_5$ has 13 pairs of Voronoi vectors pointing at 26  lattice points $xv_1+yv_2+zv_3$ with coefficients $x,y,z\in\{0,\pm 1\}$ excluding the origin $x=y=z=0$.
In any obtuse superbase $\{u_0,u_1,u_2,u_3\}$ of $\La$, all four vectors cannot have even sums of coordinates in the basis $v_1,v_2,v_3$, otherwise they cannot express the vector $v_1=(1,0,0)$.
Then at least one vector (say) $u_0$ has an odd sum $\pm 1$ or $\pm 3$.
\medskip

\textbf{The first case} is an odd sum $[u_0]=\pm 3$.
Applying reflections in the $x,y,z$-axes, we can assume that $u_0=(-1,-1,-1)=-v_1-v_2-v_3$.
Then each $u_i$, $i=1,2,3$ has no coordinate $-1$, the sum $u_0+u_i$ has coordinate $-2$, which contradicts Lemma~\ref{lem:partial_sums} saying that all partials sum of an obtuse superbase are Voronoi vectors with coordinates $x,y,z\in\{0,\pm 1\}$.
If all $u_1,u_2,u_3$ are even-sum vectors, the only remaining choice (up to permutation) is $u_1=(1,1,0)$, $u_2=(1,0,1)$, $u_3=(0,0,1)$, but all these vectors pairwisely have acute angles.
Hence one vector (say) $u_1$ has $[u_1]=1$ and we can assume that $u_0=v_1=(1,0,0)$ up to permutation. 
The sum $[u_0]+[u_1]=-3+1=-2$ can be neutralised only by Voronoi vectors $u_2,u_3$ with non-negative coordinates and $[u_2]=1=[u_3]$, so the only choice (up to a swap) is $u_2=v_2=(0,1,0)$ and $u_3=(0,0,1)$.
By reflections, this superbase $\{v_0,v_1,v_2,v_3\}$ generates all eight odd superbases in (\ref{lem:V5superbases}o), which are isometric to each other and have the coform with zeros in the top row and $p_{0i}=-v_0\cdot v_i=|v_i|^2$.
\medskip

\textbf{The second case} is an odd sum $[u_0]=\pm 1$.
Permutions and reflections in the $x,y,z$-axes allow us to assume that $u_0=(1,0,0)$.
Since $u_0$ has non-acute angles with each $u_i$, $i=1,2,3$, the first coordinates of $u_i$ is $0$ or $(-1)$. 
Since the $x$-coordinates of $u_i$ cannot have opposite signs, the projections of $u_1,u_2,u_3$ to the $(y,z)$-plane cannot have larger pairwise scalar products than the original vectors.
Hence these projections $u'_1,u'_2,u'_3$ form an obtuse superbase for the rectangular lattice $\La_2$ that they generate.
In this 2-dimensional case, all four obtuse superbases of $\La_2$ have the form $u'_2=\pm v_2$, $u'_3=\pm v_3$, $u'_1=-u'_2-u'_3$.
\medskip

Without loss of generality, assume that $u'_2=(1,0)$, $u'_3=(0,1)$.
To lift these projections to $\R^3$, if we complement both $u'_2,u'_3$ by the first coordinate $0$, we get one of the odd superbases above.
If we complement $u'_2$ by $(-1)$, then $u_2=(-1,1,0)$.
Since $u_3$ cannot have the first coordinate $(-1)$, the only choice is
$u_3=(0,0,1)$, then $u_1=-u_0-u_2-u_3=(0,-1,-1)$.
The resulting obtuse superbase $\{v_1,-v_2-v_3,v_2-v_1,v_3\}$
is $\{v_i,v_j,v_k-v_i,-v_k-v_j\}$ for $v_i=v_1$, $v_j=v_3$, $v_k=v_2$.
For any fixed $k=1,2,3$, we can choose three signs of pairwisely orthogonal basis vectors $\pm v_i,\pm v_j,\pm v_k$ in 8 ways.
Then (\ref{lem:V5superbases}e) has $3\times 8$ even obtuse superbases.
\medskip

Up to index-permutations, 
these 24 even obtuse superbases have coforms in (\ref{lem:V5superbases}o) computed from $u_0=v_i$, $u_1=v_j$, $u_2=v_k-v_i$, $u_3=-v_k-v_j$ as follows:
$p_{23}=(v_k-v_i)\cdot(v_k+v_j)=|v_k|^2$,
$p_{13}=v_j\cdot(v_k+v_j)=|v_j|^2$,
$p_{12}=v_j\cdot(v_k-v_i)=0$,
$p_{01}=-v_i\cdot v_j=0$,
$p_{02}=-v_i\cdot(v_k-v_i)=|v_i|^2$,
$p_{03}=v_i\cdot(v_k+v_j)=0$.
\medskip

The resulting coform $\mat{|v_k|^2}{|v_j|^2}{0}{0}{|v_i|^2}{0}$ can be re-written as  $\mat{0}{0}{|v_i|^2}{0}{|v_k|^2}{|v_j|^2}$ using the index-permutation induced by the composition of $0\lra 1$ 
and $0\lra 3$.
The index-permutation induced by the composition $0\lra 2$, $1\lra 3$
swaps $|v_i|^2,|v_j|^2$.
\medskip

The 24 even superbases split into three isometry classes, each having its own squared lengths $|v_i|^2,|v_j|^2,|v_i|^2+|v_k|^2,|v_j|^2+|v_k|^2$.
These unordered quadruples differ for $k=1,2,3$ if $|v_1|^2,|v_2|^2,|v_3|^2$ and their pairwise sums are all different.
\ws
\end{proof}


In 1934 the book \cite[Fig.~64 on page 170]{delone1934mathematical} gave numbers of isometry classes of obtuse superbases for Voronoi types $V_2,V_3,V_4,V_5$ as $2,3,3,1+2$, respectively, without proof.
The comment on the same page added that, for Voronoi type $V_5$, ``one class has eight superbases, each of the other two classes has three pairs of opposite superbases (six in each class)''.
In 1975 the survey \cite[Fig.~13 on page 101]{delone1975bravais} repeated the same picture with 1+2 classes for the 5th Voronoi (Dirichlet) type but added that ``there are twelve pairs of such quadrilaterals [obtuse superbases], of which the first four can differ from the second four and the third four''.
In 2009 the book \cite[p.~77]{galiulin2009crystallographic} mentioned 32 pairs of centrally symmetric obtuse superbases for a cuboid, which actually has 16 such pairs, see Example~\ref{tab:4non-isometric} below.
\medskip

Lemma~\ref{lem:V5superbases} corrects the above numbers to $1+3$ classes, where each of the three even classes in (\ref{lem:V5superbases}e) consists of eight isometric superbases, see Table~\ref{tab:4non-isometric}.
\medskip


\begin{table}[h]
\label{tab:4non-isometric}
\caption{For a primitive orthorhombic lattice, $(1+3)\times 8$ obtuse superbases split into $1+3$ isometry classes from Lemma~\ref{lem:V5superbases} and can be distinguished by lengths of vectors.
}
\medskip

\hspace*{-4mm}
\begin{tabular}{r|r|r|r}      
8 odd superbases  & 
1st even class in (\ref{lem:V5superbases}e) & 
2nd even class  in (\ref{lem:V5superbases}e) & 
3rd even class  in (\ref{lem:V5superbases}e) \\
\hline

$v_1=(\pm 1,0,0)$ & 
$v_1=(\pm 1,0,0)$ &
$v_1=(\pm 1,0,0)$ &
$v_2=(0,\pm 2,0)$ \\


$v_2=(0,\pm 2,0)$ & 
$v_2=(0,\pm 2,0)$ &
$v_3=(0,0,\pm 3)$ &
$v_3=(0,0,\pm 3)$ \\


$v_3=(0,0,\pm 3)$ &
$v_3-v_1=(\mp 1,0,\pm 3)$ & 
$v_2-v_1=(\mp 1,\pm 2,0)$ &
$v_1-v_2=(\pm 1,\mp 2,0)$ \\


$v_0=(\mp 1,\mp 2,\mp 3)$ &
$-v_3-v_2=(0,\mp 2,\mp 3)$ & 
$-v_2-v_3=(0,\mp 2,\mp 3)$ & 
$-v_1-v_3=(\mp 1,0,\mp 3)$ \\

lengths $1,2,3,\sqrt{14}$ &
lengths $1,2,\sqrt{10},\sqrt{13}$ &
lengths $1,3,\sqrt{5},\sqrt{13}$ &
lengths $2,3,\sqrt{5},\sqrt{10}$ \\

$\CF=\mat{0}{0}{0}{1}{4}{9}$ &
$\CF=\mat{0}{0}{1}{0}{9}{4}$ &
$\CF=\mat{0}{0}{1}{0}{4}{9}$ &
$\CF=\mat{0}{0}{4}{0}{1}{9}$ 
\end{tabular}
\end{table}

Lemma~\ref{lem:V5superbases} can be considered as a limit case of both Lemmas~\ref{lem:V3superbases}--\ref{lem:V4superbases}.
Indeed, 
\smallskip

\noindent
$
\{-v_1-v_2, v_1+v_3, v_2, -v_3\} 
=\{v_i,v_j,v_k-v_i,-v_k-v_j\}$, $v_i=v_2$, $v_j=-v_3$, $v_k=v_1$.

\noindent
$
\{v_0,-v_1,v_2+v_1,v_3+v_1\}
=\{v_i,v_j,v_k-v_i,-v_k-v_j\}$, $v_i=v_0$, $v_j=-v_1$, $v_k=v_2$.

\noindent
$
\{v_0+v_3, v_2+v_3, v_1, -v_3\}
=\{v_i,v_j,v_k-v_i,-v_k-v_j\}$, $v_i=v_1$, $v_j=-v_3$, $v_k=-v_0$.

\section{A root forms and a unique root invariant of a 3-dimensional lattice}
\label{sec:invariants3d}

Lemmas~\ref{lem:V1superbases}-\ref{lem:V5superbases} showed that coforms of any lattice $\La\subset\R^3$ should be considered up to different permutations for five Voronoi types. 
To reduce the ambiguity of coforms, Definition~\ref{dfn:RI} introduces below a root form $\RF(B)$ and root invariant $\RI(B)$, which will be proved to be a complete invariant of $\La\subset\R^3$ up to isometry.
\medskip

Since any obtuse superbase $B$ has only non-negative conorms, the \emph{root products} $r_{ij}=\sqrt{p_{ij}}$ are well-defined for all distinct indices $i,j\in\{0,1,2,3\}$ and have the same units as coordinates of basis vectors, for example Angstroms: $1\AA=10^{-10}$m.
The six root products can combined into a $2\times 3$ matrix called a root form $\RF=\mat{r_{23}}{r_{23}}{r_{12}}{r_{01}}{r_{02}}{r_{03}}$, which will be considered up to permutations from Lemmas~\ref{lem:V1superbases}-\ref{lem:V5superbases}.
The root invariant $\RI(B)$ will finally reduce the ambiguity $\RF(B)$ to 6, 5, 4, 4, 3 root products for Voronoi types $V_1$, $V_2$, $V_3$, $V_4$, $V_5$, respectively. 
Theorem~\ref{thm:RFinvariant} will show that the invariant $\RI(\La)$ depends only on the isometry class of $\La$.

\begin{dfn}[root form $\RF(B)$, root invariant $\RI(B)$]
\label{dfn:RI}
$\pmb{(V_5)}$
By Lemma~\ref{lem:V5superbases} any obtuse superbase $B$ of a lattice $\La\subset\R^3$ of Voronoi type $V_5$ has exactly three non-zero root products.
Up to 24 index-permutations, the \emph{root form} is $\RF(B)=\mat{0}{0}{0}{\rda}{\rdb}{\rdc}$ for any odd superbase $B$ and $\RF(B)=\mat{0}{0}{\rda}{0}{\rdb}{\rdc}$ for any even superbase $B$, where all non-zero root products are freely permutable. 
The \emph{root invariant} $\RI(B)$ is an ordered triple of the non-zero root products $\rda,\rdb,\rdc$.
\medskip

\noindent
$\pmb{(V_4)}$
For any lattice $\La\subset\R^3$ of Voronoi type $V_4$, any obtuse superbase $B$ has two zero root products in different columns.
A \emph{root form} is $\RF(B)=\mat{0}{0}{\rab}{\rda}{\rdb}{\rdc}$, where $r_{23}=0=r_{13}$, and the root products $\rab,\rda,\rdb$ are freely permutable.
The \emph{root invariant} $\RI(B)=\{(\rab,\rda,\rdb),\rdc\}$ consists of $3+1$ root products, where the triple $(\rab,\rda,\rdb)$ should be written in increasing order.
\medskip

\noindent
$\pmb{(V_3)}$
For any lattice $\La\subset\R^3$ of Voronoi type $V_3$, any obtuse superbase $B$ of $\La$ has exactly two zero root products in the same column.
A \emph{root form} is $\RF(B)=\mat{0}{\rac}{\rab}{0}{\rdb}{\rdc}$ with $r_{23}=0=r_{03}$, and $\rac,\rab,\rdb,\rdc$  are freely permutable.
The \emph{root invariant} $\RI(B)$
consists of the four non-zero root products in increasing order.
\medskip

\noindent
$\pmb{(V_2)}$
For any lattice $\La\subset\R^3$ of Voronoi type $V_2$, any obtuse superbase $B$ of $\La$ has exactly one zero root product.
A root form is $\RF(B)=\mat{0}{\rac}{\rab}{\rda}{\rdb}{\rdc}$, where
$r_{23}=0$ and the $2\times 2$ submatrix $\matfour{\rac}{\rab}{\rdb}{\rdc}$ can be changed by the symmetry group $D_4$, which
can guarantee (without changing indices for simplicity) that $\rac=\min\{\rac,\rab,\rdb,\rdc\}$ and also $\rab\leq\rdb$.
The \emph{root invariant} consists of $1+3+1$ root products: 
$\RI(B)=\{\rda,(\rac,\rab,\rdb),\rdc\}$, where $\rdc\geq\rac\leq\rab\leq\rdb$.
\medskip

\noindent
$\pmb{(V_1)}$
For any obtuse superbase $B$ of a lattice $\La\subset\R^3$ of Voronoi type $V_1$, a \emph{root form} $\RF(B)$ is the matrix $\mat{r_{23}}{r_{13}}{r_{12}}{r_{01}}{r_{02}}{r_{03}}$, where root products can be rearranged by the 24 index-permutations from Definition~\ref{dfn:index-permutations}.
A permutation of indices 1, 2, 3 as in (\ref{dfn:index-permutations}a) allows us to arrange the three columns in any order.
The composition of transpositions $0\lra i$ and $j\lra k$ for distinct $i,j,k\neq 0$ vertically swaps the root products in columns $j$ and $k$, for example apply the transposition $2\lra 3$ to the result of $0\lra 1$ in (\ref{dfn:index-permutations}b).
So we can put $r_{min}=\min\{r_{ij}\}$ into the top left position ($r_{23}$).
Then we consider the four root products in columns 2 and 3.
Keeping column 1 fixed, we can put the minimum of these four into the top middle position ($r_{13}$).
Then the resulting root products in the top row should be in increasing order.
\medskip

If the top left and top middle root products are equal ($r_{23}=r_{13}$), we can put their counterparts ($r_{01}$ and $r_{02}$) in the bottom row of columns 1,2 in increasing order.  
If the top middle and top right root products are equal ($r_{13}=r_{12}$), we can put their counterparts ($r_{02}$ and $r_{03}$) in the bottom row of columns 2 and 3 in increasing order.  
The resulting uniquely ordered matrix is the \emph{root invariant} $\RI(B)$ and can be  visualised as in Fig.~\ref{fig:forms3d}~(right) with root products instead of conorms.
\bs
\end{dfn}


\begin{lem}[equivalence of $\VF,\CF,\RF$]
\label{lem:forms3d_equiv}
For any obtuse superbase $B$, its voform $\VF(B)$, coform $\CF(B)$, and $\RI(B)$ are reconstructable from each other.
\bt
\end{lem}
\begin{proof}
The six conorms $p_{ij}$ are uniquely expressed via the seven vonorms $v_i^2,v_{ij}^2$ by formulae (\ref{dfn:forms3d}ab) and vice versa.
The root invariant $\RI(B)$ is uniquely defined by a tailored ordering of root products $r_{ij}=\sqrt{p_{ij}}$ in Definition~\ref{dfn:RI}.
\end{proof}

Important Lemmas~\ref{lem:V1superbases}--\ref{lem:V5superbases} imply in Theorem~\ref{thm:RFinvariant} below that $\RI(B)$, which was initially defined for an obtuse superbase $B$, is an isometry invariant of $\La$.
\medskip

The 2-dimensional analogue was the much simpler result in \cite[Theorem~3.7]{kurlin2022mathematics} saying that all obtuse superbases of any lattice $\La\subset\R^2$ are isometric to each other.

\begin{thm}[isometry invariance of $\RI(\La)$]
\label{thm:RFinvariant} 
If obtuse superbases $B,B'$ generate isometric lattices $\La,\La'\subset\R^3$, respectively, then $\RI(B)=\RI(B')$.
Hence $\RI$ is an isometry invariant of a lattice $\La$ and can be denoted by $\RI(\La)$.
\bs
\end{thm}
\begin{proof}
Any isometry $f$ between given lattices $\La,\La'$ maps $B$ to a new obtuse superbase $f(B)$ of $\La'$ and preserves all lengths and scalar products of vectors, so $\RF(B)=\RF(f(B))$, hence $\RI(B)=\RI(f(B))$.
Now the lattice $\La'$ has two obtuse superbases $B'$ and $f(B)$.
Lemmas~\ref{lem:V1superbases}--\ref{lem:V5superbases} explicitly described all potentially non-isometric superbases of the same lattice for five types of Voronoi domains.
\medskip

In all cases, Definition~\ref{dfn:RI} introduced the root invariant $\RI(B)$ whose root products are uniquely ordered, resolving the ambiguity of obtuse superbases.
Hence $\RI(B)=\RI(f(B))=\RI(B')$, so $\RI(\La)$ is an isometry invariant of the lattice $\La$.
\ws
\end{proof}


\begin{exa}[root invariants of orthorhombic lattices]
\label{exa:orthorhombic_lattices}
$(\pmb{oP})$
The primitive orthorhombic lattice $\La$ with edge-lengths $0\leq a\leq b\leq c$ has the obtuse superbase $v_1=(a,0,0)$, $v_2=(0,b,0)$, $v_3=(0,0,c)$, $v_0=(-a,-b,-c)$, whose root form is $\RF(\La)=\mat{0}{0}{0}{a}{b}{c}$, so the root invariant is $\RI(\La)=(a,b,c)$.
If we re-order vectors, columns of $\RF(\La)$ are re-ordered accordingly, but $\RI(\La)$ remains the same.
Another obtuse superbase $v_1=(a,0,0)$, $v_2=(0,b,0)$, $v'_3=(-a,0,c)$, $v'_0=(0,-b,-c)$ has $\RF(\La)=\mat{0}{0}{a}{0}{b}{c}$, but $\RI(\La)=(a,b,c)$ is the same.
\medskip

\noindent
$(\pmb{tP})$ 
For a primitive tetragonal lattice, set $a=b$ in the case above.
\medskip

Let all orthorhombic lattices below have a base cube with
sides $2a\leq 2b\leq 2c$. 
\medskip

\noindent
$(\pmb{oS})$
A base-centred orthorhombic lattice $\La$ has the obtuse superbase $v_1=(2a,0,0)$, $v_2=(-a,b,0)$, $v_3=(0,0,c)$, $v_0=(-a,-b,-c)$, whose root form is $\RF(\La)=\mat{0}{0}{a\sqrt{2}}{a\sqrt{2}}{\sqrt{b^2-a^2}}{c}$.
The root invariant is $\RI(\La)=\{(a\sqrt{2},a\sqrt{2},\sqrt{b^2-a^2}),c\}$ where the $\sqrt{b^2-a^2}$ should move to the first place if $a\sqrt{2}>\sqrt{b^2-a^2}$, $b<a\sqrt{3}$.
\medskip

\noindent
$(\pmb{oF})$
A face-centred orthorhombic lattice $\La$ has the obtuse superbase $v_1=(a,b,0)$, $v_2=(a,-b,0)$, $v_3=(-a,0,c)$, $v_0=(-a,0,-c)$, whose root form is $\RF(\La)=\mat{\sqrt{b^2-a^2}}{a}{a}{\sqrt{c^2-a^2}}{a}{a}$.
If $b<a\sqrt{2}$, the root invariant is $\RI(\La)=\RF(\La)$, otherwise $\RI(\La)$ is obtained from $\RF(\La)$ by swapping the first column with the last column.
\medskip

\noindent
$(\pmb{oI})$
For a body-centred orthorhombic lattice $\La$, assume that $a^2+b^2\geq c^2$.
Then $\La$ has the obtuse superbase $v_1=(a,b,-c)$, $v_2=(a,-b,c)$, $v_3=(-a,b,c)$, $v_0=(-a,-b,-c)$, and $\RF(\La)=\mat{\sqrt{a^2+b^2-c^2}}{\sqrt{a^2-b^2+c^2}}{\sqrt{-a^2+b^2+c^2}}{\sqrt{a^2+b^2-c^2}}{\sqrt{a^2-b^2+c^2}}{\sqrt{-a^2+b^2+c^2}}$.
Due to $a\leq b\leq c$, the root products are increasing in each row, so $\RI(\La)=\RF(\La)$.
\medskip

\noindent
$(\pmb{tI})$ 
For a body-centred tetragonal lattice, set $a=b$ in the case above.
\bs
\end{exa}

\section{Root invariants classify all 3-dimensional lattices up to isometry}
\label{sec:classification3d}

Proposition~\ref{prop:superbases_lattices} substantially reduces the ambiguity of  lattice representations by their bases.
Any fixed lattice $\La\subset\R^3$ has infinitely many (super)bases but only a few obtuse superbases, maximum 32 (or two non-isometric classes) 
in Lemma~\ref{lem:V5superbases}.

\begin{prop}[obtuse superbases of isometric lattices in $\R^3$]
\label{prop:superbases_lattices} 
Lattices in $\R^3$ are isometric if and only if any of their obtuse superbases $B,B'$ are isometric to each other or to a couple of obtuse superbases in one of Lemmas~\ref{lem:V1superbases}--\ref{lem:V5superbases}.
\bt
\end{prop}
\begin{proof}
Part \emph{only if} ($\Rightarrow$): any isometry $f$ between lattices $\La,\La'$ maps any obtuse superbase $B$ of $\La$ to the obtuse superbase $f(B)$ of $\La'$.
Then $B'$ and $f(B)$ are isometric to each other or to obtuse superbases listed in one of Lemmas~\ref{lem:V1superbases}--\ref{lem:V5superbases} for the Voronoi type of the given isometric lattices $\La\cong\La'$.  
\medskip

Part \emph{if} ($\Leftarrow$): the given conditions on $B,B'$ mean that there is an isometry $B\to B'$ extending to an isometry of their lattices $\La\to\La'$, or $B,B'$ are isometric to a couple of obtuse superbases in one of Lemmas~\ref{lem:V1superbases}--\ref{lem:V5superbases}, so $\La\cong\La'$.
\ws
\end{proof}

Proposition~\ref{prop:superbases_lattices} above formalises the key difference between dimensions 2 and 3 for an isometry classification of lattices.
The 2D analogue in \cite[Theorem~3.7]{kurlin2022mathematics} says that any lattices $\La\subset\R^2$ are isometric if and only if any their superbases are isometric.
Proposition~\ref{prop:superbases_lattices} needs much more sophisticated  Lemmas~\ref{lem:V1superbases}--\ref{lem:V5superbases}, because lattices in $\R^3$ can have several non-isometric superbases.

\begin{lem}[superbase reconstruction]
\label{lem:superbase_reconstruction}
An obtuse superbase of any lattice $\La\subset\R^3$ can be reconstructed up to isometry from its root invariant $\RI(\La)$.
\bt
\end{lem}
\begin{proof}
The root invariant $\RI(\La)$ can be lifted to a $2\times 3$ matrix of a root form $\RF(\La)$ for each of five Voronoi types of lattices in Definition~\ref{dfn:RI}.
The positions of root products $r_{ij}=\sqrt{-v_i\cdot v_j}$ in $\RF(\La)$ allow us to compute the lengths $|v_i|$ from the formulae of Definition~\ref{dfn:forms3d}, for example $|v_0|=\sqrt{r_{01}^2+r_{02}^2+r_{03}^2}$.
Up to rigid motion in $\R^3$, one can fix $v_0$ along the positive $x$-axis in $\R^3$.
\medskip

The angle $\angle(v_i,v_j)=\arccos\dfrac{v_i\cdot v_j}{|v_i|\cdot|v_j|}\in[0,\pi)$ between the vectors $v_i,v_j$ can be found from the vonorms $v_i^2,v_j^2$ and root product $r_{ij}=\sqrt{-v_i\cdot v_j}$.
The found length $|v_1|$ and angle $\angle(v_0,v_1)$ allow us to fix $v_1$ in the $xy$-plane of $\R^3$.
The vector $v_2$ with the known length $|v_2|$ and two angles $\angle(v_0,v_2)$ and $\angle(v_1,v_2)$ has two positions that are isometric by the mirror reflection in the $xy$-plane.
\ws
\end{proof}

\begin{thm}[3D lattices/isometry $\lra$ root invariants]
\label{thm:classification3d}
Any lattices $\La,\La'\subset\R^3$ are isometric if and only if their root invariants coincide: $\RI(\La)=\RI(\La')$.
\bt
\end{thm}
\begin{proof}
The part \emph{only if} ($\Rightarrow$) is Theorem~\ref{thm:RFinvariant} implying that any isometric lattices $\La,\La'$ have equal root invariants: $\RI(\La)=\RI(\La')$.
The part \emph{if} ($\Leftarrow$) follows from Lemma~\ref{lem:superbase_reconstruction} by reconstructing a superbase of $\La$ from its root invariant $\RI(\La)$. 
\ws
\end{proof}

\begin{cor}[3D lattices/similarity $\lra$ proportional $\RI$]
\label{cor:similarity3d}
Any lattices 
in $\R^3$ are related by similarity (a composition of isometry and uniform scaling) if and only if their root invariants are proportional by a factor $s>0$. 
\bt
\end{cor}
\begin{proof}
Scaling a lattice $\La\subset\R^3$ by a factor $s>0$ multiplies all root products $r_{ij}$, hence all components of $\RI(\La)$, by $s$.
The corollary follows from Theorem~\ref{thm:classification3d}.
\ws
\end{proof}

\begin{exa}[non-isometric lattices with $DC^7=0$]
\label{exa:dc7=0}
Fig.~\ref{fig:DC7examples} shows that we cannot freely permute  vonorms or conorms (equivalently, root products) without changing the isometry class of a lattice.
The voforms in Fig.~\ref{fig:DC7examples} differ by a single transposition $10\lra 12$ for the vonorms $v_{12}^2=v_{03}^2$ and $v_{23}^2=v_{01}^2$.
This transposition is not among the 24 index-permutations from Definition~\ref{dfn:forms3d}.
The coforms in Fig.~\ref{fig:DC7examples} are computed from the voforms by formulae~(\ref{dfn:forms3d}b).
These coforms include different conorms, for example value 5 appear in $\CF(\La)$ but not in $\CF(\ti\La)$.
Then $\RI(\La)\neq\RI(\La')$ define non-isometric lattices $\La\not\cong\ti\La$ by Theorem~\ref{thm:classification3d}. 
\medskip

In these lattices $\La,\ti\La\subset\R^3$ the origin $0$ has the same distances $|v_0|$, $|v_1|$, $|v_2|$, $|v_3|$, $|v_{12}|$, $|v_{23}|$, $|v_{13}|$ to its seven closest Voronoi neighbours.  
Hence the function $DC^7$ taking the Euclidean distance between these 7-dimensional distance vectors \cite{andrews2019space} vanishes
for $\La,\ti\La$.
Our colleagues Larry Andrews and Herbert Bernstein quickly checked that $\La,\ti\La$ can be distinguished by the 8th distance from the origin to its 8th closest neighbour.   
However, the example in Fig.~\ref{fig:DC7examples} can be extended to an infinite 6-parameter family of non-isometric lattices $\La,\ti\La$ with $DC^7(\La,\ti\La)=0$ as follows.
\medskip

\begin{figure}[h]
\label{fig:DC7examples}
\caption{The lattices $\La,\ti\La$ defined by the coforms $\CF(\La),\CF(\ti\La)$ are not isometric due to $\RI(\La)\neq\RI(\ti\La)$ but the origin $0$ has the same distances to its seven closest neighbours in both $\La,\ti\La$.}
\includegraphics[width=\textwidth]{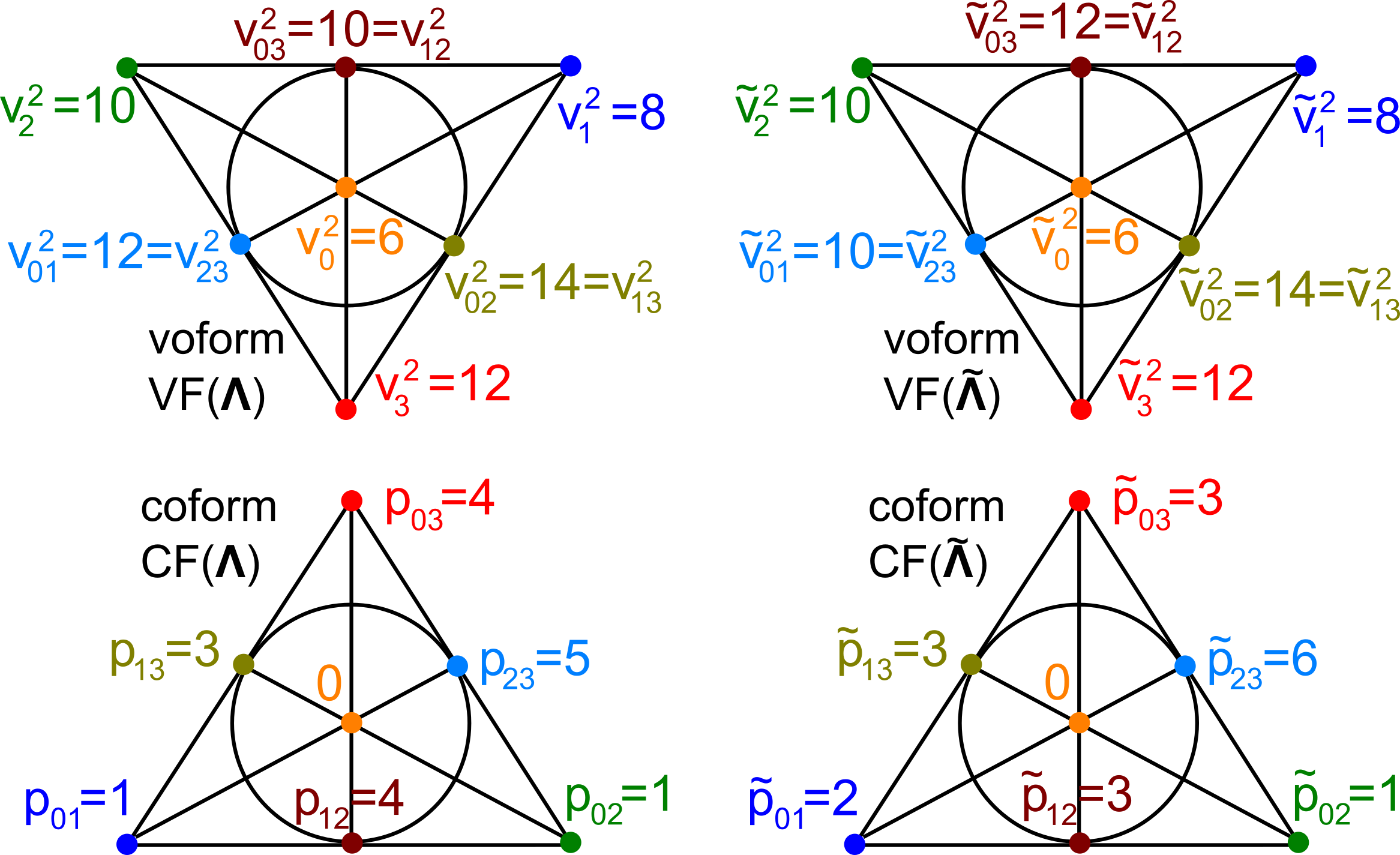}
\end{figure} 

Add any conorms $q_{ij}\geq 0$ to $\CF(\La),\CF(\ti\La)$ in the `conorm-wise' way.
Formulae~(\ref{dfn:forms3d}a) imply that the voforms $\VF(\La),\VF(\ti\La)$ consist of the same 7 numbers:
$$\La:\left\{\begin{array}{ll} 
v_0^2 & =(p_{01}+q_{01})+(p_{02}+q_{02})+(p_{03}+q_{03})
=1+4+1+q_{01}+q_{02}+q_{03},\\
v_1^2 &=(p_{01}+q_{01})+(p_{12}+q_{12})+(p_{13}+q_{13})=1+4+3+q_{01}+q_{12}+q_{13},\\
v_2^2 &=(p_{02}+q_{02})+(p_{12}+q_{12})+(p_{23}+q_{23})=1+4+5+q_{02}+q_{12}+q_{23},\\
v_3^2 &=(p_{03}+q_{03})+(p_{13}+q_{13})+(p_{23}+q_{23})=4+3+5+q_{03}+q_{13}+q_{23},\\
v_{01}^2 &=(p_{02}+q_{02})+(p_{03}+q_{03})+(p_{12}+q_{12})+(p_{13}+q_{13})= \\
& =1+4+4+3+q_{02}+q_{03}+q_{12}+q_{13}=
12+(q_{02}+q_{03}+q_{12}+q_{13}), \\
v_{02}^2 &=(p_{01}+q_{01})+(p_{03}+q_{03})+(p_{12}+q_{12})+(p_{23}+q_{23})= \\
& =1+4+4+5+q_{01}+q_{03}+q_{12}+q_{23}=
14+(q_{01}+q_{03}+q_{12}+q_{23}), \\
v_{03}^2 &=(p_{01}+q_{01})+(p_{02}+q_{02})+(p_{13}+q_{13})+(p_{23}+q_{23})= \\
& =1+1+3+5+q_{01}+q_{02}+q_{13}+q_{23}=
10+(q_{01}+q_{02}+q_{13}+q_{23}); \\
\end{array}\right.$$
$$\ti\La:\left\{\begin{array}{ll}  
\ti v_0^2 & =(\ti p_{01}+q_{01})+(\ti p_{02}+q_{02})+(\ti p_{03}+q_{03})=2+1+3+q_{01}+q_{02}+q_{03},\\
\ti v_1^2 & =(\ti p_{01}+q_{01})+(\ti p_{12}+q_{12})+(\ti p_{13}+q_{13})=2+3+3+q_{01}+q_{12}+q_{13},\\
\ti v_2^2 & =(\ti p_{02}+q_{02})+(\ti p_{12}+q_{12})+(\ti p_{23}+q_{23})=1+3+6+q_{02}+q_{12}+q_{23},\\
\ti v_3^2 & =(\ti p_{03}+q_{03})+(\ti p_{13}+q_{13})+(\ti p_{23}+q_{23})=3+3+6+q_{03}+q_{13}+q_{23}, \\
\ti v_{01}^2 &=(\ti p_{02}+q_{02})+(\ti p_{03}+q_{03})+(\ti p_{12}+q_{12})+(\ti p_{13}+q_{13})= \\
& =1+3+3+3+q_{02}+q_{03}+q_{12}+q_{13}=
10+(q_{02}+q_{03}+q_{12}+q_{13}),\\
\ti v_{02}^2 &=(\ti p_{01}+q_{01})+(\ti p_{03}+q_{03})+(\ti p_{12}+q_{12})+(\ti p_{23}+q_{23})= \\
& =2+3+3+6+q_{01}+q_{03}+q_{12}+q_{23}=
14+(q_{01}+q_{03}+q_{12}+q_{23}), \\
\ti v_{03}^2 &=(\ti p_{01}+q_{01})+(\ti p_{02}+q_{02})+(\ti p_{13}+q_{13})+(\ti p_{23}+q_{23})= \\
& =2+1+3+6+q_{01}+q_{02}+q_{13}+q_{23}=
12+(q_{01}+q_{02}+q_{13}+q_{23}). \\
\end{array}\right.$$
Notice that almost all vonorms coincide: $v_i^2=\ti v_i^2$ and $v_{02}^2=\ti v_{02}^2$ except the couple of swapped values:
$v_{01}^2=\ti v_{03}^2$ and $v_{03}^2=\ti v_{01}^2$.
So both lattices $\La,\ti\La$ have the same ordered distances from the origin to its seven closest neighbours: $DC^7(\La,\ti\La)=0$.
\medskip

Now we show that the new coforms $\CF(\La),\CF(\ti\La)$ lead to different root invariants for almost all free parameters $q_{ij}\geq 0$ in the generic case of Lemma~\ref{lem:V1superbases} when all conorms are positive.
Under $4!=24$ index-permutations from Definition~\ref{dfn:forms3d}, any two conorms from a common column remain in together in a (possibly another) column.
The coforms $\CF(\La),\CF(\ti\La)$ have the following column sums
$$\La:\left\{\begin{array}{l} 
p_{23}+p_{01}=(5+q_{23})+(1+q_{01})=6+q_{23}+q_{01},\\
p_{13}+p_{02}=(3+q_{13})+(1+q_{02})=4+q_{13}+q_{02},\\
p_{12}+p_{03}=(4+q_{12})+(4+q_{03})=8+q_{12}+q_{03};
\end{array}\right.$$
$$\ti\La:\left\{\begin{array}{l} 
\ti p_{23}+\ti p_{01}=(6+q_{23})+(2+q_{01})=8+q_{23}+q_{01},\\
\ti p_{13}+\ti p_{02}=(3+q_{13})+(1+q_{02})=4+q_{13}+q_{02},\\
\ti p_{12}+\ti p_{03}=(3+q_{12})+(3+q_{03})=6+q_{12}+q_{03}.
\end{array}\right.$$

Two sums from the above triples coincide: $p_{13}+p_{02}=\ti p_{13}+\ti p_{02}$ for any $q_{ij}\geq 0$.
Since the first sums and third sums clearly differ, the above triples of sums can coincide only if the remaining pairs of sums are swapped, so $p_{23}+p_{01}=\ti p_{12}+\ti p_{03}$ and $p_{12}+p_{03}=\ti p_{12}+\ti p_{03}$, which both are equivalent to $q_{23}+q_{01}=q_{12}+q_{03}$.
If $q_{23}+q_{01}\neq q_{12}+q_{03}$, the above triples of sums differ, so $\CF(\La),\CF(\ti\La)$ are not related by index-permutations.
The underlying lattices $\La,\ti\La$ are not isometric by Theorem~\ref{thm:classification3d}.
To distinguish the lattices $\La\not\cong\ti\La$ in this 6-parameter family by 8 or more distances from the origin to its neighbours, a theoretical proof is needed.
\bs
\end{exa}

Lemma~\ref{lem:bounds_products} below implies that the root products $r_{ij}$ continuously change under perturbations of an obtuse superbase measured in the Minkowski metric $M_{\infty}$.

\begin{lem}[bounds for root products {\cite[Lemma~7.3]{kurlin2022mathematics}}]
\label{lem:bounds_products}
Let vectors \\ $u_1,u_2,v_1,v_2\in\R^n$ have a maximum length $l$, have non-positive scalar products $u_1\cdot u_2,v_1\cdot v_2\leq 0$, and 
$|u_i-v_i|\leq\de$ for $i=1,2$.
Then 
$$|u_1\cdot u_2-v_1\cdot v_2|\leq 2l\de \qquad\text{ and }\qquad
|\sqrt{-u_1\cdot u_2}-\sqrt{-v_1\cdot v_2}|\leq\sqrt{2l\de}.
\eqno{\blacktriangle}$$
\end{lem}

The next paper \cite{kurlin2022easily} will prove a stronger continuity result by defining metrics on root invariants.
Justifying metric axioms will be much harder than in $\R^2$ \cite[section~5]{kurlin2022mathematics}, because we need to glue five Voronoi type subspaces of $\LIS(\R^3)$ in a non-trivial way not covered by the classical theory \cite[Part I, Lemma 5.24]{bridson2013metric}. 
These continuous metrics will define real-valued chiralities of 3D lattices by continuously measuring a deviation from a higher-symmetry neighbour as in \cite[section~6]{kurlin2022mathematics}.
\medskip

The more recent Pointwise Distance Distributions \cite{widdowson2021pointwise} are continuous, complete for distance-generic crystals and helped establish 
the \emph{Crystal Isometry Principle} saying that all real periodic crystals can be distinguished up to isometry by their geometric structures of atomic centres without chemical data. 
Hence all periodic crystals live in the common Crystal Isometry Space (CRISP), which can be projected to the Lattice Isometry Space $\LIS(\R^3)$ parameterised in Problem~\ref{pro:map}.
\medskip

The companion papers in dimension 2 \cite{bright2021geographic} and 3 \cite{bright2021welcome} discuss many continuous maps of real crystal lattices from the Cambridge Structural Database. 
\medskip

Many thanks to all colleagues who read early drafts for their valuable time. 

\begin{acknowledgements}
This research was supported by the £3.5M EPSRC grant `Application-driven Topological Data Analysis' (2018-2023), the £10M Leverhulme Research Centre for Functional Materials Design (2016-2026) and the Royal Academy of Engineering Fellowship `Data Science for Next Generation Engineering of Solid Crystalline Materials' (2021-2023).
\end{acknowledgements}


\bibliographystyle{spmpsci}      
\bibliography{lattices3Dmaths}   


\renewcommand{\thesection}{\Alph{section}}
\setcounter{section}{0}

\section{Proof of reduction: any 3D lattice has an obtuse superbase}
\label{sec:proofs}

The previous paper \cite[Appendix~A]{kurlin2022mathematics} includes  some basic definitions and proofs of past results outlined by Delone, Conway and Sloane \cite{conway1992low}.
This appendix corrects (in the next update) the example in \cite[Fig.~8 in section~7]{conway1992low} used in the proof of \cite[Theorem~8]{conway1992low}.
Below we give a more detailed argument for this Reduction Theorem~\ref{thm:reduction} by using Lemma~\ref{lem:reduction} as a typical reduction step.

\begin{lem}[reduction]
\label{lem:reduction}
Let $B=(v_0,v_1,v_2,v_3)$ be any superbase of a lattice $\La\subset\R^3$.
For any distinct $i,j,k,l\in\{0,1,2,3\}$, let the new superbase vectors be 
$u_i=-v_i$, $u_j=v_j$, $u_k=v_{ik}=v_i+v_k$, $u_l=v_{il}=v_i+v_l$.
Then all vonorms remain the same or swap their places, and the only change is $u_{ij}^2=v_{ij}^2-4\ep$, where $\ep=v_i\cdot v_j$.
The conorms $q_{\bullet}$ of the new vectors $u_{\bullet}$ are updated 
as in Fig.~\ref{fig:forms3d_reduction}
$$q_{ij}=\ep,\; 
q_{jk}=p_{jk}-\ep,\; 
q_{jl}=p_{jl}-\ep,\; 
q_{ik}=p_{il}-\ep,\; 
q_{il}=p_{ik}-\ep,\; 
q_{kl}=p_{kl}+\ep. 
\leqno{(\ref{lem:reduction})}
$$
\end{lem}

\begin{proof}
If initial vectors $v_{\bullet}$ form a superbase, which means that $v_i+v_j+v_k+v_l=0$, then so do the new vectors: $u_i+u_j+u_k+u_l=(-v_i)+v_j+(v_i+v_k)+(v_i+v_l)=0$.
\medskip

\begin{figure}[h]
\caption{Lemma~\ref{lem:reduction} for $i=1$, $k=2$, $j=3$, $l=0$ says that the new superbase $u_1=-v_1$, $u_2=v_{12}$, $u_3=v_3$, $u_0=v_{01}$ has the new voform $\VF$ and coform $\CF$ shown above.}
\includegraphics[width=\textwidth]{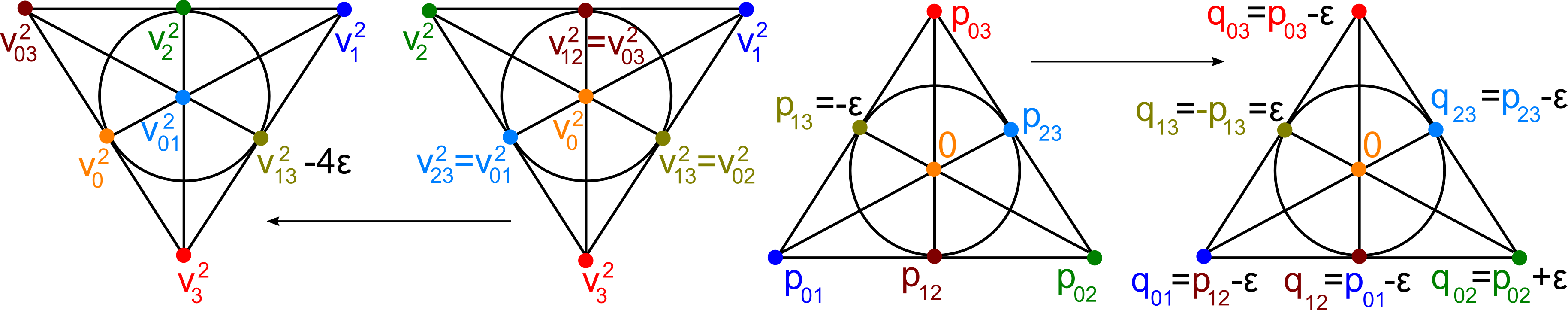}
\label{fig:forms3d_reduction}
\end{figure}

For the new superbase $u_i=-v_i$, $u_j=v_j$, $u_k=v_{ik}$, $u_l=v_{il}$,
two vonorms remain the same: $u_i^2=v_i^2$ and $u_j^2=v_j^2$.
Two pairs of vonorms swap values: $u_k^2=v_{ik}^2$, $u_{jl}^2=u_{ik}^2=(u_i+u_k)^2=v_k^2$ and $u_l^2=v_{il}^2$, $u_{jk}^2=u_{il}^2=(u_i+u_l)^2=v_l^2$.
The final vonorm is
$$u_{ij}^2=u_{kl}^2=(v_j-v_i)^2=(v_i+v_j)^2-4v_i\cdot v_j=v_{ij}^2+4p_{ij}=v_{ij}^2-4\ep, \text{ see Fig.~\ref{fig:forms3d_reduction}.}$$
We similarly check (\ref{lem:reduction}) illustrated in Fig.~\ref{fig:forms3d_reduction} for $i=1$, $k=2$, $j=3$, $l=0$.

\noindent
$q_{ij}=-u_i\cdot u_j=v_i\cdot v_j=-p_{ij}=\ep$

\noindent
$q_{jk}=-u_j\cdot u_k=-v_j\cdot(v_i+v_k)=-v_i\cdot v_j-v_j\cdot v_k=p_{jk}-\ep$

\noindent
$q_{jl}=-u_j\cdot u_l=-v_j\cdot(v_i+v_l)=-v_i\cdot v_j-v_j\cdot v_l=p_{jl}-\ep$

\noindent
$q_{ik}=-u_i\cdot u_k=v_i\cdot(v_i+v_k)=v_i\cdot(-v_j-v_l)=-v_i\cdot v_l-v_i\cdot v_j=p_{il}-\ep$

\noindent
$q_{il}=-u_i\cdot u_l=v_i\cdot(v_i+v_l)=v_i\cdot(-v_j-v_k)=-v_i\cdot v_k-v_i\cdot v_j=p_{ik}-\ep$

\noindent
$q_{kl}=-u_k\cdot u_l=(v_i+v_k)\cdot(v_i+v_l)=v_i\cdot(-v_i-v_j-v_k)-v_k\cdot v_l
=v_i\cdot v_j+p_{kl}
=p_{kl}+\ep$.
\medskip

Notice that the conorm $p_0$ at the centre of $\CF$ remains zero by formula~(\ref{lem:vonorms<->conorms}b):
$$4p_{0}=u_i^2+u_j^2+u_k^2+u_l^2-u_{ij}^2-u_{ik}^2-u_{il}^2
=v_i^2+v_j^2+(v_i+v_k)^2+(v_i+v_l)^2-(v_j-v_i)^2-v_k^2-v_l^2=$$
$$=v_i^2+v_j^2+(v_i^2+2v_iv_k+v_k^2)+(v_j+v_k)^2-(v_i^2-2v_iv_j+v_j^2)-v_k^2-(v_i+v_j+v_k)^2=$$
$$=v_i^2+v_j^2+v_k^2+2v_iv_j+2v_jv_k+2v_iv_k-(v_i+v_j+v_k)^2
=0.$$
Hence all central conorms $p_0$ in \cite[Fig.~5]{conway1992low} should be 0.
\ws
\end{proof}

\begin{proof}[of Theorem~\ref{thm:reduction} for $n=3$]
We will reduce any superbase $B=(v_0,v_1,v_2,v_3)$ of a lattice $\La\subset\R^3$ to make all conorms $p_{ij}$ non-negative.
Starting from a largest negative conorm $p_{ij}=-\ep<0$, change the superbase by Lemma~\ref{lem:reduction}.
This reduction leads to the positive conorm $q_{ij}=\ep$, not zero as in \cite[Fig.~4(b)]{conway1992low}.
\medskip

Four other conorms decrease by $\ep>0$ and can potentially become negative, which requires a new reduction by Lemma~\ref{lem:reduction} and so on.
To prove that the reduction process always finishes, notice that six vonorms keep or swap their values, but one vonorm always decreases by $4\ep>0$.
Every reduction can make superbase vectors only shorter, but not shorter than a minimum distance between points of $\La$.
The angle between $v_i,v_j$ can have only finitely many values when lengths of $v_i,v_j$ are bounded.
Hence the scalar product $\ep=v_i\cdot v_j>0$ cannot converge to 0.
Since every reduction makes one partial sum $v_S$ shorter by a positive constant, while other six vectors $v_S$ keep or swap their lengths, the reductions by Lemma~\ref{lem:reduction} should finish in finitely many steps.
\ws
\end{proof}

A reduction of lattice bases for real crystals has many efficient implementations.
Theoretical estimates for reduction steps are discussed in \cite{nguyen2009low}.

\end{document}